\input amstex.tex
\input amsppt.sty   
\magnification 1200
\hsize=6.2 true in
\vsize = 8.6 true in
\nologo
\NoRunningHeads        
\parskip=\medskipamount
        \lineskip=2pt\baselineskip=18pt\lineskiplimit=0pt
       
        \TagsOnRight
        \NoBlackBoxes

         \topmatter
        \title
       Quasi-periodic solutions to a nonlinear Klein-Gordon equation\\
       with a decaying nonlinear term
      \endtitle
\author
         W.-M.~Wang        \endauthor        
\address
{CNRS and Department of Mathematics, Cergy Paris Universit\'e, 95302 Cergy-Pontoise Cedex, France}
\endaddress
 \email
{wei-min.wang\@math.cnrs.fr}
\endemail
\abstract
We present a set of time quasi-periodic solutions to a nonlinear Klein-Gordon equation with a decaying nonlinear term
on the torus in arbitrary dimensions. This generalizes the bifurcation method  developed in \cite{W2}.


\endabstract

 \bigskip\bigskip
        \bigskip
        \toc
        \bigskip
        \bigskip 
        \widestnumber\head {Table of Contents}
        \head 1. Introduction and statement of the Theorem
        \endhead
        \head 2. The good linear solutions
        \endhead
        \head 3. Extraction of parameters -- solving the $\quad\quad\quad$ 
         $Q$-equations
        \endhead
        \head 4. The first step -- solving the $P$-equations
        \endhead
        \head 5. The $\theta$ estimates
        \endhead
         \head 6. Proof of the Theorem
        \endhead
         \endtoc
        \endtopmatter
        \vfill\eject
        \bigskip
\document
 
     
\head{\bf 1. Introduction and statement of the Theorem}\endhead
We consider {\it real} valued solutions to a nonlinear Klein-Gordon equation (NLKG) on the $d$-torus $\Bbb T^d=[0, 2\pi)^d$:
$$
\frac{\partial^2u}{\partial t^2}-\Delta u+u+(\Cal M*u)^{p+1}*\Cal M=0,\tag 1.1
$$
where $p\in\Bbb N$, is arbitrary; considered as functions on $\Bbb R^d$, $u$ and $\Cal M$ satisfy: 
$u(\cdot, x)=u(\cdot, x+2j\pi)$ and $\Cal M(x)=\Cal M(x+2j\pi)$ for all $j\in\Bbb Z^d$. The role of $\Cal M$
is to regularize certain range of frequencies.

We start from the linear, second order in time, equation: 
$$
\frac{\partial^2u}{\partial t^2}-\Delta u+u=0,\tag 1.2
$$
and define the wave operator $D$: 
$$D:=\sqrt{-\Delta+1}.\tag 1.3$$ 
Using the Fourier series, it follows that the spectrum of $D$: 
$$\sigma (D)=\{\sqrt{|j|^2+1}\,|\,j\in\Bbb Z^d\}.$$
In $d\geq 2$, the spectrum is degenerate and the gap between non-equal eigenvalues shrinks to zero.


Denote $|j|^2$ by $j^2$, the solutions to (1.2) are linear combinations of cosine and sine functions of the form:
$$\cos ({-(\sqrt{j^2+1})t+j\cdot x})\tag 1.4$$ and 
$$\sin ({-(\sqrt{j^2+1})t+j\cdot x}),\tag 1.5$$
where 
$\cdot$ is the usual scalar product. These solutions are, generally speaking, quasi-periodic in time (``periodic" with several frequencies). 

After the addition of the nonlinear terms, it is natural to investigate the
bifurcation of these quasi-periodic solutions. We use the space-time approach initiated in \cite{W2}, cf. also \cite{W3}
for a review of this approach.

\subheading {1.1 Statement of the Theorem}
To streamline the presentation, we seek solutions to (1.1) which have space-time reflection
symmetry: $u(t, x)=u(-t, -x)$. Therefore, we assume $\Cal M (x)=\Cal M (-x)$. 

Let $u^{(0)}$ be an even solution of $b$ frequencies, to the linear equation (1.2): 
$$u^{(0)}(t, x)=\sum_{k=1}^b a_k\cos (-(\sqrt{j_k^2+1})t+j_k\cdot x).\tag 1.6$$
For the nonlinear construction, it is useful to add a dimension for each frequency in time and view $u^{(0)}$
as a function on $\Bbb T^b\times\Bbb T^d : =\Bbb T^{b+d}\supset\Bbb T^d$. 
Let $$\omega^{(0)}:=\{\sqrt{j_k^2+1}\}_{k=1}^{b},\, (j_k\neq 0),$$
be the $b$-dimensional frequency vector. 
Henceforth $u^{(0)}$ may be written in the form: 
$$
\aligned u^{(0)}(t, x)&=\sum_{k=1}^b a_k\cos (-(\sqrt{j_k^2+1})t+j_k\cdot x)\\
:&=\sum_{k=1}^b\hat u^{(0)}(\mp e_{k}, \pm j_k)\cos (\mp e_k\cdot\omega^{(0)}t\pm j_k\cdot x), \endaligned
$$
where $e_{k}= (0, 0, ...1, .., 0)\in\Bbb Z^b$ is a unit vector, with the only non-zero component in the $k$th direction, 
and $$\hat u(-e_{k}, j_k)=\hat u(e_{k}, -j_k)=a_k/2.$$
We say that ${\hat u}^{(0)}$ has support
$$\text{supp } {{\hat u}^{(0)}}=\{(\mp e_{k}, \pm j_k), k=1,...,b\}\subset\Bbb Z^{b+d},\tag 1.7$$
where $j_k\neq j_{k'}$ if $k\neq k'$.

For the nonlinear equation (1.1), we seek quasi-periodic solutions with $b$ frequencies in the form of a space-time cosine series:
$$
u(t, x)=\sum_{(n,j)\in\Bbb Z^{b+d}}\hat u(n, j)\cos ({n\cdot\omega t}+j\cdot x), \tag 1.8
$$
satisfying $\hat u(n, j)=\hat u(-n, -j)$ and with the frequency $\omega\in\Bbb R^b$ to be determined. 
We note that the corresponding linear solution $u^{(0)}$
has {\it fixed} frequency $$\omega=\omega^{(0)}=\{\sqrt {j_k^2+1}\}_{k=1}^b\in\Bbb R^b,$$ 
which are eigenvalues of the wave operator $D$ defined in (1.3). 

The cosine space is an ``invariant" subspace for the NLKG in (1.1).
Restricting to the cosine series amounts to restricting to the cosine-cosine sector and economizes considerably the notations.

We say that a solution to the linear equation (1.2) is {\it good} if its spatial frequencies satisfy the non-degeneracy conditions (i-iii) in sect.~2.1.
Considering $\{j_1, j_2, ..., j_b\}$ as a point in $(\Bbb Z^d)^b$, it suffices to say here that there are {\it infinite} number of good frequencies, 
and in fact the set of good frequencies have {\it positive density} in $(\Bbb Z^d)^b$.
(Cf. the Lemma in sect.~2, and Remark. 2 above sect.~2.2.) 

Let $\Vert u^{(0)}\Vert =\delta$, and $\hat {\Cal M}$ satisfy: 
$$\aligned \hat {\Cal M}(m)&=1, \qquad\qquad\qquad\, |m|\leq e^{|\log\delta|^{1/2}};\\
&\leq e^{-|m|},\quad\quad\qquad\, \text{ } |m|>e^{|\log\delta|^{1/2}}.\endaligned \tag $\flat$ $$
Below is our main result: 
\proclaim
{Theorem} 
Let $$u^{(0)}(t, x)=\sum_{k=1}^b a_k\cos (-(\sqrt{j_k^2+1})t+j_k\cdot x),$$ 
be a good solution to the linear equation (1.2), satisfying the non-degeneracy conditions (i-iii), $a=\{a_k\}_{k=1}^b\in (-\delta,\delta)^b\backslash \{0\}=\Cal B(0, \delta)$ and $p$ even. 
Assume $b>C_p d$, where $C_p$ only depends on $p$. 
Then for all $\epsilon\in (0, 1)$,  there exists $\delta_0>0$, such that 
for all $\delta\in (0, \delta_0)$,  there is a 
Cantor set $\Cal G\subset \Cal B(0, \delta)$ with 
$$\text{meas }\Cal G/(2\delta)^b\geq 1-\epsilon,$$ 
and a diffeomorphism: $a\mapsto\omega (a)$ on $\Cal B(0, \delta)$.
For all $a\in\Cal G$, there is an analytic quasi-periodic solution to the nonlinear equation (1.1) of the form (1.8),
with $\omega=\omega(a)=\{\omega_k(a)\}^b_{k=1}$ satisfying
$$\omega_k=\sqrt{j_k^2+1}+\Cal O(\delta^{p}),\, k=1, 2, ..., b,$$
and $$|u(t, x)-\sum_{k=1}^b a_k\cos (-\omega_kt+j_k\cdot x)|=\Cal O(\delta^{p}),$$
for all $(t, x)\in\Bbb R\times \Bbb T^b$. 
 \endproclaim
 \noindent{\it Remark.}  
The concept of good linear solutions remains valid for odd $p$. 
The additional assumption of even $p$ is to ensure 
amplitude-frequency modulation at the leading order $\Cal O(\delta^{p})$. This is a sufficient but {\it not} necessary condition. 
The condition of large $b$, namely $b>C_pd$, is imposed in order that certain determinants are not identically zero, cf. Proof of Lemma~4.2. This is the same reason as in \cite{W2}. 
 \subheading {1.2 About the characteristic and $\Cal M$} 
Using (1.8) in (1.2), yields the support of $\hat u$; setting $\omega=\omega^{(0)}$ leads to:
$$\aligned \Cal C:=&\{(n,j)\in{\Bbb Z}^{b+d}|-(n\cdot\omega^{(0)} )^2+j^2+1=0\}\\
  =&\{(n,j)\in\Bbb Z^{b+d}|\pm n\cdot\omega^{(0)} +\sqrt{j^2+1}=0\}\\
:= &\Cal C^{+}\cup\Cal C^-.\endaligned\tag 1.9$$ 
The set $\Cal C$ is an infinite set of resonances, and is one of the main difficulties here.
We call $\Cal C$, the characteristic, and see it as the restriction to $\Bbb Z^{b+d}$ of the corresponding  
hyperboloids on $\Bbb R^{b+d}$. 

We divide the Fourier space $\Bbb Z^{b+d}$ according to the scales $N=|(n, j)|$:
\item{(i)} $\log N\ll |\log\delta|$;
\item{(ii)}   $\log N\sim |\log\delta|$;
\item{(iii)}   $\log N\gg|\log\delta|$.

For the small scales (i), using algebraic method and number theory, we establish good separation property (by the name of connected sets), in Proposition~2.1 and Lemma~5.2, on the characteristic $\Cal C$.
We note that the regularization does not affect the nonlinear term at these scales,  since $\hat {\Cal M}=1$ from ($\flat$).
For large scales (iii), the separation lemma, Lemma~20.14 \cite {B4} near $\Cal C$, suffices (after appropriate
rescaling by $\delta$). The intermediate scales (ii) are a new phenomenon. Lemma~20.14 \cite {B4}, which needs parameters $\omega$ 
of order $1$, is not applicable here, since $\omega^{(0)}$ are {\it fixed}.  The generalization of Lemma~20.14 \cite {B4} to arbitrary dimensions
for fixed frequencies currently remains open. 
\footnote{In two dimensions, however, this can be generalized \cite{W4}, and should
lead to the existence of quasi-periodic solutions for the (usual) polynomial nonlinearity 
$u^{p+1}$, instead of the $(\Cal M* u)^{p+1}*\Cal M$ in (1.1).}
On the other hand, the proof of separation property at small scales is rather general, applicable to complicated spectra
and is valid at arbitrary dimensions. So it could be of independent interest, and we present it here by proving the Theorem.


We could regularize the problem, using a convolution potential which only acts on the scales (ii). However, since the main novelty 
of the paper is the analysis of small scales, for simplicity, we took the $\Cal M$ in ($\flat$) instead, which is smoothing also for large scales (iii).  
Fast polynomial decay suffices, but the exponential decay in ($\flat$) does facilitate the expositions in sects. ~5.6 and 6. 
\subheading{1.3 Some background} 
Quasi-periodic solutions have been previously constructed in one dimension 
with positive mass $m$. In that case, the linear Klein-Gordon equation:
$$
\frac{\partial^2u}{\partial t^2}-\frac{\partial^2u}{\partial x^2}+mu=0,
$$ 
gives rise to an eigenvalue set $\{\sqrt {j^2+m}, j\in\Bbb Z \}$ close to the set of integers, see \cite{B2}\cite{P}, and in a related context \cite{CY}\cite{K1}\cite{Way}.
For almost all $m$, this 
set is rationally independent. This property does not have higher dimensional analogues and is a serious obstacle. 
(In the special time periodic case, solutions have been constructed in higher dimensions in \cite{B1}.)  For NLKG with a multiplicative, non-constant potential,
see \cite{BeBo}. For quasi-periodic solutions to NLS, see \cite{W2} \cite{PP}, see also the related works \cite{B3,~4}\cite{EK}; for parameter dependent nonlinear beam equations, see \cite{GY}\cite{EGK}. 
\subheading {1.4 Comparison of NLKG with NLS }
The NLKG in (1.1) gives rise to an infinite dimensional dynamical system. 
Studying such a system usually requires certain separation property in
order to approximate it by ``direct sums" of finite dimensional systems.  
In fact, {\it separation} on the characteristic is indispensable to
the existence of KAM-type solutions for nonlinear PDEs -- without it, there could be Arnold diffusion like scenario, cf. \cite{CKSTT}
for related finite time results for NLS.
Furthermore, this separation property is used to  {\it extract parameters} from the nonlinear term for later analysis.  

The characteristic in (1.9) defines hyperboloids; while for NLS
the bi-characteristics are paraboloids, which are limit-elliptic, cf. (1.5, 1.6) in \cite{W2}.  Due to the 
convexity or ellipticity of the Laplacian, 
separation property for NLS is attained 
by considering intersections of hyperplanes, and
occasionally ellipsoids, which are compact \cite{W2}.

By contrast, NLKG is {\it hyperbolic}, separation entails intersections of quadratic hyper-surfaces,
which are generically non-compact, cf. sects.~5.1-5.4. It is much more difficult to discern that a large system of quadratic polynomial equations
in many variables has no solution. 
To circumvent this difficulty, a pre-selection of 
the initial frequencies $\omega^{(0)}=\{\omega_k^{(0)}\}_{k=1}^b$ is made along number theoretical considerations. More precisely, the 
frequencies are chosen to be square roots of distinct {\it square-free} integers. This is the key new feature,
compared to NLS, cf. proof of the Lemma in sect. 2. 

As a direct consequence of this frequency selection, there is the linear independence, cf. e.g., \cite{Ro} for a proof (for notational simplicity, $\omega$ is generally written for $\omega^{(0)}$) 
:$$\Vert n\cdot\omega\Vert_{\Bbb T}\neq 0,$$
where $n\neq 0$ and $\Vert\, \Vert_{\Bbb T}$ denotes the distance to the integers; as well as the {\it quadratic} non-equality:
$$\Vert \sum_{k<\ell}n_{k}n_{\ell}\omega_k\omega_\ell \Vert_{\Bbb T}\neq 0,\tag 1.10$$
where $\sum_{k<\ell} |n_{k}n_{\ell} | \neq 0$. (See the elementary derivation in sect.~2.1.)
The latter implies that:  
$$\Vert \sum_{k,\ell} n_{k}n_{ \ell} \omega_k\omega_{\ell}\Vert _{\Bbb T}=0\, \Longleftrightarrow \sum_{k, \ell} n_{k}n_{ \ell} \omega_k\omega_{\ell}=n^2_{r} \omega^2_r,\tag 1.11$$
for some $r\in\{1,2, ..., b\}$.

The linear independence is the usual one; the quadratic non-equality is new and takes care of hyperbolicity,
moreover it 
{\it doubles} as a small-divisor lower bound:
If $$\pm n\cdot\omega+\sqrt{j^2+1}\neq 0,\tag 1.12$$
then 
$$|\pm n\cdot\omega+\sqrt{j^2+1}|>c|n|^{-q},\, n\neq 0,\tag 1.13$$
where $c$, $q>0$, by using the simultaneous Diophantine approximation result in \cite{Schm},
cf. the very beginning of sect. 4.1.
We note that these are new types of small-divisors, which do not appear in NLS. 
\subheading {1.5 Effective resultant analysis}
Realizing the aforementioned finite dimensional approximations, 
leads to the concept of {\it good} linear solutions $u^{(0)}$. The {\it effective resultant} analysis introduced in \cite{W2} is used 
to control the size of the finite dimensional systems.  

For NLS, the effective resultants $\Cal D$ can be reached with relative ease, since it mostly involves intersections of hyper-planes. 
The condition $\Cal D\neq 0$ may be imposed, 
cf. the proof of Lemma~2.5 in \cite{W2}. 
For NLKG this is not feasible, particularly so in the proof of Lemma~5.2, as it involves intersections of hyperboloids. 
Properties of square roots of square-free integers are used instead to arrive at the effective resultants $\Cal D$. 
This arithmetic condition seems natural, indeed it would be difficult otherwise. 
The variety defined by $\Cal D=0$ is then analyzed with the additional help of the curvature near the origin, 
introduced by the mass term $1$ in the wave operator. 
(For more details, see the proofs of Proposition~2.1 and Lemma~5.2.) 
This generalizes the non-perturbative bifurcation analysis introduced in \cite{W2} to NLKG. 
\subheading  {1.6 The nonlinear matrix equation} 
Using the ansatz in (1.8), the NLKG in (1.1) is directly transformed into a nonlinear matrix equation in the ``Fourier coefficients" $\hat u (n, j)$.
(See (2.2), sects.~3 and 4 for more details.) The equation is divided into two parts using a Lyapunov-Schmidt $P$- and $Q$-equations 
decomposition. The domain of the $Q$-equations is the set supp ${\hat u}^{(0)}$ defined in (1.7); while that of the $P$-equations, 
the complement set. 
The $Q$-equations are used to solve for the 
frequencies $\omega=\omega(a)$, which remain real (since we work in the real) and permit amplitude-frequency modulation; 
while the $P$-equations, the Fourier coefficients $\hat u$, iteratively using a multi-scale Newton scheme. As in \cite{W2}, 
due to the resonances, the equation is linearized at the unperturbed solution $u^{(0)}$, instead of at $0$.
Ensuring invertibility of the linearized operator  then leads to the non-degeneracy conditions (i-iii) on $u^{(0)}$ in sect.~2.1. 
\subheading  {1.7 The multiscale analysis}
The iteration in sect.~5.6, relies on the invertibility of the linearized operators,
which is obtained by {\it multiscale} analysis. Here multiscale means both
in space and time (Fourier space).  From this standpoint, the earlier sections address the initial 
scale (or initial scales). The idea of multiscale analysis (MA) originated in the study of Anderson localization in \cite{FS}.  

Central is the control of resonances at each scale. More precisely, for our problem, say scale $N$ means that 
the Fourier variables are restricted to the cube $\Lambda_N=[-N, N]^{b+d}$ and the linearized operator is restricted to $\Lambda_N$. 
One covers $\Lambda_N$ by smaller cubes $\Lambda_{N'}$ with appropriately chosen $N'\ll N$. MA requires that there is 
{\it separation} of the resonant cubes $\Lambda_{N'}$, or there are only few of them. The version of MA that is relevant here is \cite{BGS} and \cite{W1} -- {\it few} means that at each scale $N$, 
there are only {\it sublinear in} $N$, resonant $N'$-cubes. For small scales $N$, there is separation, as mentioned in sect.~1.4; while for larger scales, few resonances, from Diophantine property and the decaying
nonlinear term. 

\noindent{\it Remark.} The multiscale analysis method can also prove that the spectrum of the Floquet operator is pure point, as in e.g., \cite {FSW}, \cite{BGS} and {\cite{W1}. 
\subheading {1.8 Organization of the paper}
In sect.~2, we define {\it good} linear solutions and connected sets. The non-degeneracy conditions are then used to 
bound the size of connected sets on the characteristic $\Cal C$ -- the size of the ``finite dimensional"
systems. The $Q$-equations are solved in sect.~3, 
leading to amplitude-frequency modulation. In sect.~4, the first corrections to the linear solutions are 
obtained using the Newton scheme. The linear analysis needed in all subsequent Newton iterations are done
in sect.~5. The proof concludes by constructing the quasi-perioidic solutions in sect.~6. 
\subheading {Notations}

We summarize below some of the notational conventions:

\noindent-- The dimension $d$, the degree of nonlinearity $p$ and the number of basic frequencies $b$ are fixed. The set 
$\{j_k\}_{k=1}^b$ is a fixed subset of $\Bbb Z^d$.

\noindent-- The letter $u$ denotes a function on $\Bbb T^{b+d}$, $\hat u$ its Fourier series.
The hat is generally dropped and $u$ is written for $\hat u$, which are 
functions on $\Bbb Z^{b+d}$. 

\noindent-- The letters $n$ and $\nu$ denote vectors in $\Bbb Z^b$; while $j$ and $\eta$ vectors in $\Bbb Z^d$.  

\noindent-- The dot $\cdot$ denotes the usual scalar product in Euclidean space. To simplify notations, one writes $j^2$ for $j\cdot j$ etc. 

\noindent-- The norm $\Vert\,\Vert$ stands for the $\ell^2$ or operator norm; while $|\,|$ for the sup-norm or the length of a vector in a finite 
dimensional vector space or the number of elements in a given set.

\noindent-- An identically zero function $f$ is denoted by $f \equiv  0$; the negation $f\not\equiv 0$.

\noindent-- Given two positive quantities, $A, B>0$, $A\asymp B$ signifies $cB<A<CB$ for some $0<c<C$; while $A\lesssim B$, 
$A<CB$.

\noindent-- Large positive constants are generally denoted by upper case letters such as $C$, $C'$ etc.; while small ones $c$, $c'$, $\epsilon$ etc.
Unless indicated otherwise, they are not the same and may vary from statement to statement.
\smallskip

\head{\bf 2. The good linear solutions}\endhead  
Let $*$ denote convolution on $\Bbb Z^{b+d}$: 
$$[\hat V*\hat W](x)=\sum_y\hat V(x-y)\hat W(y).\tag 2.1$$
When $V$ is an even function on $\Bbb T^{b+d}$, its Fourier series satisfies, moreover, 
$$\hat V(x-y)=\hat V(y-x),$$
so (2.1) defines a {\it self-adjoint} operator. 
Using the ansatz (1.8), (1.1) becomes
$$\text{ diag } [-(n\cdot\omega)^2+j^2+1]\hat u+(\hat\Cal M{\hat u})^{*{ (p+1) }}\hat\Cal M=0.\tag 2.2$$
From now on we work with (2.2), and for simplicity we drop the hat and write $u$ for $\hat u$ and $\Cal M$ for $\hat \Cal M$.

We seek solutions close to the linear solution $u^{(0)}$ of $b$ frequencies, 
$\text{supp } {u}^{(0)}=\{(\mp e_{k}, \pm j_k), k=1,...,b\},$ with frequencies
$\omega^{(0)}=\{\sqrt{j_k^2+1}\}_{k=1}^{b}$ ($j_k\neq 0$) 
and small amplitudes $a=\{a_k\}_{k=1}^b$ satisfying $\Vert a\Vert=\Cal O(\delta)\ll 1$.
Denote the left side of (2.2) by $F(u)$.  

Linearizing at $u^{(0)}$, we 
are led to study the linearized operator $F'(u^{(0)})$ on $\ell^2(\Bbb Z^{b+d})$ with  
$$F'=D'+A,\tag 2.3$$
where 
$$
D'=\text{ diag } [-(n\cdot\omega)^2+j^2+1]\tag 2.4$$
and
$$
A=(p+1)\Cal M[(u^{(0)})^{*p}*]\Cal M.\tag 2.5$$
Since $u^{(0)}$ is even, $A$ and hence the linearized operator $F'$ in (2.3) are self-adjoint. 
Moreover $A$ is a T\"oplitz matrix in the $n$-direction.

In order for  $u^{(0)}$ to bifurcate to a nearby solution $u$ to the nonlinear equation in (1.1), algebro-geometric
conditions will be imposed on the support of  $u^{(0)}$. This is in the spirit of \cite{W2}. The algebraic 
aspect here is, however, new, and originates from arithmetic considerations for the temporal frequencies of $u^{(0)}$.

\subheading {2.1 The good linear solutions}
To define good $u^{(0)}$, we need to analyze the
matrix $A$ defined in (2.5).  More specifically, the 
{\it structure} of the matrix, which is determined by supp $u^{(0)}$.  
($\hat \Cal M$ affects only the size of the entries.)

Let
$$\Gamma=\text{supp }[(u^{(0)})^{*p}]\backslash \{(0, 0)\}=\{(\nu, \eta)\}\subset\Bbb Z^{b+d},\tag 2.6$$ 
with $$\text{supp } u^{(0)}=\{(\mp e_k, \pm j_k)\}_{k=1}^b,$$
where $e_k\in \Bbb Z^b$ and $j_k\in\Bbb Z^d$. 
From the definition,  $(\nu, \eta)\in \Gamma$ are of the form 
$$(-\sum_{k=1}^b m_ke_k, \sum_{k=1}^b m_kj_k),$$ 
where $\sum_{k=1}^b |m_k|\leq p$. So 
$$\eta=-\sum_{i=1}^b \nu_ij_i.$$ 

More generally, for any fixed $R\in \Bbb N$, consider the set 
$$\tilde\Gamma_R
=\bigcup_{r=1}^{R} \text{supp }[(u^{(0)})^{*pr}]\backslash \{(0, 0)\},\tag$\sharp$ $$
and again use $(\nu, \eta)$ to denote an element of $\tilde\Gamma_R$: $(\nu, \eta)\in \tilde\Gamma_R$.  If
$(\nu, \eta)\in\tilde \Gamma_R$, then $(\nu, \eta)$ is of the form 
$$(-\sum_{k=1}^b m_ke_k, \sum_{k=1}^b m_kj_k),$$ 
where $\sum_{k=1}^b |m_k|\leq pR.$ So there remains the relation: 
$$\eta=-\sum_{i=1}^b \nu_ij_i, \tag $\sharp\sharp$ $$ 
for $(\nu, \eta)\in \tilde\Gamma_R$. 

The vector $\eta=\eta (j_1, j_2, ..., j_b)$ is considered as a function from 
$(\Bbb Z^d)^b$ to $\Bbb Z^d$. More precisely, for a given $\nu=\{\nu_i\}_{i=1}^b\in\Bbb Z^b$, $\eta=-\sum_{i=1}^b\nu_ij_i$ 
is a function from $(\Bbb Z^d)^b$ to $\Bbb Z^d$. 
There is the basic relation:
$$\nu=0 \Longleftrightarrow\,  \eta\equiv 0, \tag $\bigstar$ $$
which can be seen as follows. From  ($\sharp\sharp$), if $\nu=0$, then $\eta\equiv 0$; if $\eta\equiv 0$ and $\nu\neq 0$, then there exists 
$k\in\{1, 2, ...b\}$ such that $\nu_{k}\neq 0$. Set all $j_{k'}=0$ for $k'\neq k$, then 
$\eta\neq 0$ for $j_k\neq 0$, which is a contradiction. (See also Lemma 2.1 in \cite{W2}.)  
 \hfill $\square$

For the Klein-Gordon equation, because of the restriction derived in (1.11), $\eta$ which is a function of at most $2$ variables, $$\eta=\eta (j_1, j_2, ..., j_b)=\tilde \eta(j_k, j_\ell),$$ $k, \ell\in\{1, 2, ..., b\}$, plays an 
important role.
For notational simplicity, we write $\eta$ for $\tilde\eta$: $$\eta(j_k, j_\ell):=\tilde \eta(j_k, j_\ell).$$ 

Below $(\nu, \eta)$ is again considered as a point in 
$\Bbb Z^{b+d}$. Fix $R=2d+1$ in ($\sharp$). The following definition is to ensure that when {\it restricted} to 
$\tilde\Cal C:=\Cal C\backslash \text{supp } u^{(0)}$, where $\Cal C$ is the characteristic defined in (1.9) and $\text{supp } u^{(0)}$ in (1.7), the matrix $A$ defined in (2.5) is block-diagonal with blocks of size at most 
$2d\times 2d$. 

\noindent{\bf Definition.} $u^{(0)}$  a solution of $b$ frequencies $j_1$, $j_2$, ..., $j_b$, to the linear equation in (1.5) is {\it good} if the following 
three conditions are satisfied:
\item{(i)} If $b\geq d+1$, any $d$ vectors in the set $\{j_k \}_{k=1}^b$ are linearly independent. 
For all $k$, $k=1, 2, ..., b$, define the set 
$$J_k=\{j_{k'}-j_k|k'=1, ..., b, k'\neq k\},$$ 
any $d$ vectors in $J_k$ are linearly independent.
(If $b\leq d$, there is no condition (i).)
\item{(ii)} The integers $(j^2_k+1)$, $k=1$, $2$, ..., $b$, are distinct:
$$1<j_1^2+1<j_2^2+1<\cdots<j_b^2+1,$$
and square-free.
\item{(iii)} For all given $k\in\{1, 2, ..., b\}$ and $m\in\Bbb Z\cap [-p, p]\backslash\{0\}$, consider the set of $(\nu, \eta)\in\tilde\Gamma_{2d+1}$ with 
$$\nu=-me_k+m_\ell e_\ell,$$
where $\ell\in\{1, 2, ..., b\}$, $m_\ell\in\Bbb Z$, $|m_\ell|\leq 2pd$, and 
$$\aligned \eta=&mj_k-m_\ell j_\ell\\
:=&\eta_{\ell, m_\ell}\neq 0.\endaligned$$
For each $\eta$, define $L$ to be 
$$L=2m\eta\cdot j_k+(m^2-m_\ell^2):=L(\ell, m_\ell).$$
Denote by $P(\ell, m_\ell)$ the corresponding $d$-dimensional hyperplane in $\Bbb R^d$:
$$2\eta\cdot j+L=0,\tag *$$ 
where $\eta=\eta_{\ell, m_\ell}$ and $L=L(\ell, m_\ell)$. 

Let $\sigma$ be any set of $(\ell, m_\ell)$ with $2d$ elements, $\ell\in\{1, 2, ..., b\}$,  $m_\ell\in\Bbb Z$, $|m_\ell|\leq 2pd$, such that there exists $(\tilde\ell, m_{\tilde\ell} )\in\sigma$ with $m_{\tilde\ell}\neq \pm m$, then 
$$\bigcap_{\sigma} P(\ell, m_\ell)=\emptyset.$$

\noindent{\it Remark.} 
Instead of $1$, the mass may be fixed at any square-free integers-- the proof of the Theorem should be the same.
The conditions (i, iii)  are to be considered on $\Bbb R^d$. 

As mentioned in sect.~1, condition (ii) implies the usual linear independence: 
$$\Vert \sum_{k=1}^bn_k\omega_k\Vert_{\Bbb T}\neq 0,\tag \dag$$
where $n_k\in\Bbb Z$, $\sum_k|n_k|\neq 0$;
as well as the quadratic non-equality: 
$$\Vert \sum_{k, \ell; k< \ell}n_{k}n_{\ell}\omega_k\omega_\ell \Vert_{\Bbb T}\neq 0,\tag\dag\dag$$
where $\sum_{k<\ell} |n_{k}n_{\ell}|\neq 0$.

Property ($\dag$) follows from basic algebra (\cite{Ro}). 
We prove ($\dag\dag$) by contradiction. If the opposite of ($\dag\dag$)  holds, then 
$\sum_{k, \ell; k< \ell}n_{k}n_{\ell}\omega_k\omega_\ell \in\Bbb Z$, where $\sum_{k<\ell} |n_{k}n_{\ell}|\neq 0$.
So $$(\sum_k n_k\omega_k)^2=2\sum_{k, \ell; k< \ell}n_{k}n_{\ell}\omega_k\omega_\ell+\sum_k n_k^2\omega_k^2\in\Bbb Z^+\cup\{0\}.$$
Hence $\sum_k n_k\omega_k\in\sqrt{\Bbb Z^+}\cup\{0\},$ which is a contradiction, by using the $b+1$ term version of ($\dag$), if needed.
This yields ($\dag\dag$).  \hfill $\square$

The following indicates that the above two conditions are viable.

\proclaim {Lemma} There is an infinite number of $(j_1, j_2, ...,j_b)\in(\Bbb Z^{d})^b$ which satisfy the non-degeneracy conditions (i-iii).
\endproclaim
Since the proof is slightly lengthy, we first prove: 
\proclaim {Lemma'} There is an infinite number of $(j_1, j_2, ...,j_b)\in(\Bbb Z^{d})^b$ which satisfy the non-degeneracy conditions (i) and (iii).
\endproclaim
\demo{Proof} 
The first part of condition (i) is satisfied if the $d\times d$ determinant $D$ of any $d$ vectors
in the set $\{j_k\}_{k=1}^b$, satisfies $D\neq 0$. 
Likewise, for a fixed $k$, the second part is satisfied if the $d\times d$ determinant $D'$ of any $d$ vectors in the set $J_k$, satisfies $D'\neq 0$.
Since clearly $D\not\equiv 0$, $D'\not\equiv 0$, 
$D'\neq 0$ and $D''\neq 0$ define Zarisky open sets. The intersection of all such sets formed above, $\Cal D'$,  is therefore Zarisky 
open and contains an infinite number of integers. On $\Cal D'$, (i) is satisfied. 

For (iii), since $(0, 0)\notin\tilde\Gamma_{2d+1}$, $\eta\neq 0$ define sets of co-dimension $1$.   
Fixing $m$, $k$, we note that if $\ell=k$ (if $m_\ell=0$, set $\ell=k$), then 
$$P(k, m_k)\cap P(k, m'_{k})=\emptyset,\tag 2.7$$ if $m_k\neq m'_{k}$. This is because $P(k, m_k)$ is
the plane defined by the equation:
$$2\eta_{k, m_k}\cdot j +L(k, m_k)=0,$$
equivalently by $$2j_k\cdot j+2mj_k^2+(m+m_k)=0,$$
which is incompatible with the equation for $P(k, m'_{k})$ if $m_k\neq m'_{k}$.
Similarly,  when $\ell\neq k$,
$$P(\ell, m_\ell)\cap P(\ell, m'_{\ell})\cap P(\ell, m''_{\ell}) =\emptyset,\tag 2.8$$ 
if $m_\ell\neq m'_{\ell}$ and  $m'_\ell\neq m''_{\ell}$ and  $m''_\ell\neq m_{\ell}$. 

Below we may assume $b\geq d+1$, as otherwise 
(2.7) implies that $k$, $j_k$ may appear at most once and (2.8) implies that
for each $\ell$, $\ell\neq k$, $j_\ell$ may appear at most twice in order to have 
non-empty intersection. So
$$\bigcap_{\sigma} P(\ell, m_\ell)=\emptyset,$$
if $|\sigma|=2d$ and $b\leq d$. 
Let $\sigma'$ be a subset of $\sigma$ with $d$ elements: 
$$\sigma'\subset \sigma\backslash\{(\tilde \ell, m_{\tilde\ell})\}, \, |\sigma'|=d.$$
The corresponding set of $\eta$ is then:  
$$\{\eta_{\ell, m_\ell}=mj_k-m_\ell j_\ell; (\ell, m_\ell)\in\sigma'\}.$$
Using the same argument as above, one may assume that 
the above set of $\eta$ depends on at least $d$ variables in 
$\{j_\kappa\}_{\kappa=1}^b$, as otherwise if $\sigma$ does not contain such a subset, then
$$\bigcap_{\sigma\backslash\{(\tilde\ell, m_{\tilde\ell} )\}}P(\ell, m_\ell)=\emptyset$$
from (2.7, 2.8).

Let 
$$\tilde \sigma=\sigma'\cup\{(\tilde\ell, m_{\tilde\ell})\}.$$
One stipulates that the $(d+1)\times (d+1)$ determinant 
of  the set of vectors $$\{(2\eta_{\ell, m_\ell}, L(\ell, m_\ell)); (\ell, m_\ell)\in\tilde\sigma\},$$
$$\det (2\eta, L)\not\equiv 0. \tag 2.9$$There are two cases. 

If $\tilde \ell=k$ or if $m_{\tilde\ell}=0$, set $j_k=0$.
The first part of condition (i) gives the linear independence of 
$$\{\eta_{\ell, m_\ell}; (\ell, m_\ell)\in\sigma'\}=\{-m_\ell j_\ell; (\ell, m_\ell)\in\sigma'\}.$$
Since
$$L(\tilde\ell, m_{\tilde\ell})\neq 0,$$ this yields
$$\det (2\eta, L)\neq 0.$$  
If $\tilde \ell\neq k$, set $\eta_{\tilde\ell}=0$. 
Since $L(\tilde\ell, m_{\tilde\ell})\neq 0$,  either using the first part of condition (i) directly or after setting 
$j_k=0$ (when all $\eta_{\ell, m_\ell}$, $(\ell, m_\ell)\in\sigma'$ depend on two variables), we obtain  
$$\det (2\eta, L)\neq 0.$$ 
So $$\bar D=\det (2\eta, L)\not\equiv 0.$$ 
Let $D(k, m, \tilde\sigma)$ be the set in $(\Bbb R^d)^b$ defined by $\bar D\neq 0$. Then 
$D(k, m, \tilde\sigma)$ is Zarisky open containing an infinite number of integers. 

More generally, for any fixed $k$, $m$, let $\tilde\sigma$ be a set of $(d+1)$ 
elements $(\ell_i, m_{\ell_i})$, $i=1, 2, ..., d+1$, such that the corresponding 
subset of $\eta$ of $d$ elements: 
$$\eta_{\ell_i, m_{\ell_i}}=mj_k-m_{\ell_i}j_{\ell_i},$$ $i=1, 2, ..., d$,
depends on at least $d$ variables in $\{j_\kappa\}_{\kappa=1}^b$
and such that $m_{\ell_{d+1}}\neq \pm m$. 
Then the $(d+1)\times (d+1)$ determinant 
$$\bar D=\det (2\eta, L)\not\equiv 0$$ by using the same argument. 
Let $D(k, m, \tilde\sigma)$ be the set in $(\Bbb R^d)^b$ defined by $\bar D\neq 0$. Then 
the set 
$$\Cal D'':=\bigcap_{k, m, \tilde\sigma}D(k, m, \tilde\sigma)\neq\emptyset,$$
is Zarisky open and contains an infinite number of integers. On $\Cal D''$,
(iii) is verified.
The set $$\Cal D=\Cal D'\cap \Cal D''$$ is therefore Zarisky 
open containing an infinite number of integers. On $\Cal D$, (i) and (iii) are satisfied. 
\hfill$\square$
\enddemo
We are left to deal with (ii). It is well-known that there is an infinite number of square-free integers. 
In fact, for $N\gg 1$, in the set $\{1, 2, ... ,N\}$, there are $6N/\pi^2+\Cal O(\sqrt N)$ square-free integers. 
One of the complications here is that we are working with a set of {\it sums of squares} of integers, and {\it not} the set of integers itself.  
To complete the proof of the Lemma, we need the following:
\proclaim{Lemma 0} Denote the set of square-free integers by $S_{qf}$. Assume $0\neq m\in S_{qf}$, then 
$$\{n^2+m; n\in\Bbb Z\}\cap S_{qf}$$ is an infinite set; moreover
$$|\{n^2+m; n\in\Bbb Z\}\cap S_{qf}\cap \{1, 2, ..., N\}|=\Cal O(\sqrt{N}), \, N\gg 1.$$
\endproclaim
\demo{Proof}  Since $m\in S_{qf}$, $$\text{gcd}\{n^2+m; n\in\Bbb Z\}\in S_{qf}, $$ 
and $f(n)=n^2+m\in \Bbb Z[n]$ has 
no repeated roots for $m\neq 0$, it follows from \cite{N} that $\{n^2+m; n\in\Bbb Z\}\cap S_{qf}$ is an infinite set. That
it has positive density follows from \cite{Es}. 
\hfill$\square$
\enddemo

\demo{Proof of Lemma}
Call an integer vector $V$ in $\Bbb Z^b$ square-free if all its components are square-free. 
We are left to show that there is an infinite number of 
$$(j_1, j_2, ..., j_b)\in (\Bbb Z^d)^b\cap \Cal D, \tag $\diamondsuit$ $$
satisfying $$ 1<j_1^2+1<j_2^2+1< ... < j_b^2+1, \tag $\diamondsuit\diamondsuit$ $$ and
such that the integer vector
$$\Bbb Z^b\ni V=(j_1^2+1, j_2^2+1, ..., j_b^2+1)$$
is square-free. 
Define the set $S_d$: 
$$\Bbb Z^d\supset S_d=\{(x_1, x_2, ..., x_d)\in\Bbb Z^d|\sum_{i=1}^d x_i^2+1\in S_{qf}\},$$
and $S_d^b\subset (\Bbb Z^d)^b$ the product set.
Let $\pi_i$ be the projection onto the $i^{\text {th}}$ copy of $\Bbb Z^d$, $i=1, 2, ..., b$. It suffices that $$\pi_i \{S_d^b\cap \Cal D\}  \tag $\diamondsuit\diamondsuit\diamondsuit$ $$
is an infinite set for all $i=1, 2, ..., b$. Below we describe such a selection process. 

The set $\Cal D=\Cal D'\cap\Cal D''$;
$\Cal D'$ can be described by polynomial non-equalities of degrees at most $d$; while 
the set $\Cal D''$ at most degree $d+2$. 
These polynomials are in $bd$ variables: 
$$j_{1,1}, j_{1, 2}, ..., j_{1, d}, ..., j_{k, i}, ...,  j_{b, d}; k=1, 2, ..., b, i=1, 2, ..., d,$$ 
where $j_{k, i}$ denotes the $i$th component of $j_k\in\Bbb Z^d\subset\Bbb R^d$. There are finite number
(depending only on $d, p, b$), $N(d, p, b)$ of such polynomials. Denote this set of polynomials by $\Cal P$. 

We prove by induction. First consider $j_1$,
in the order $$j_{1,1}, j_{1, 2}, ..., j_{1, d}.$$ 
For simplicity of notation, set $$x_i=j_{1, i}, i=1, 2, ..., d.$$ 
Assume that the variable $x_1$ appears in $N_1$ polynomials, $P_1, P_2, ..., P_{N_1}$ of degrees at most $d+2$. Consider $P_1$:
$$0\not\equiv P_1\in \Bbb R[X_1, X_2, ..., X_d; Y],$$
where $Y$ denotes the variables $j_2, ..., j_b$.  
Consider $P_1$ as an element of $\Bbb R[X_2, \dots X_n; Y][X_1]$. It has at most $d+2$ roots in $\Bbb R[X_2, \dots, X_n; Y]$,
hence at most $d+2$ roots in $\Bbb R$. So the set 
$$\Sigma'_{1,1} = \{x_1 \in \Bbb R, P_1(x_1, X_2, \dots X_n; Y) \equiv 0\}$$ has at most $d+2$
elements. Clearly same consideration holds for $P_2, ..., P_{N_1}$. Denote by $\Cal P_1$ the set of polynomials $P_1, P_2, ..., P_{N_1}$,  we have
$$\Sigma_{1,1} = \{x_1 \in \Bbb R, P_i(x_1, X_2, \dots X_n; Y) \equiv 0, \text { for some }P_i\in\Cal P_1, i=1, 2, ..., N_1\}$$ has at most $(d+2)N_1\leq (d+2)N(p, q, d)$ elements.
Clearly $$\Sigma_{1, 2}(x_1) , ..., \Sigma_{1, d}(x_1, x_2, ...,  x_{d-1})$$ can be constructed similarly for $$x_1\notin\Sigma_{1, 1}, x_2\notin\Sigma_{1, 2}(x_1), ..., x_{d-1}\notin\Sigma_{1, d-1}(x_1, x_2, ..., x_{d-2}),$$
and $$|\Sigma_{1, 2}(x_1)|, ...., |\Sigma_{1, d}(x_{d-1})|\leq (d+2)N(p, q, d),$$ 
{\it uniformly} in  $$x_1\notin\Sigma_{1, 1}, x_2\notin\Sigma_{1, 2}(x_1), ..., x_{d-1}\notin\Sigma_{1, d-1}(x_1, x_2, ..., x_{d-2}).$$

This construction extends to $j_2$, ..., $j_b$ and produces the sets $$\Sigma_{k, \ell}(\{j_{k', \ell'}; k'\leq k, \ell'<\ell\})$$ for $k=2, ..., b$ and $\ell=1, 2, ..., d$.
Call $\Sigma_{k, \ell}$ the non-admissible sets.  Then 
$$|\Sigma_{k, \ell}(\{j_{k', \ell'}; k'\leq k, \ell'<\ell\}) |\leq  (d+2)N(p, q, d),$$
for all $k=1, 2, ..., b$ and $\ell=1, 2, ..., d$,
when the arguments are {\it not} in the non-admissible sets. 
Clearly, if fixing a point $(j_1, j_2, ..., j_b)\in(\Bbb Z^d)^b$ in the order 
$$j_{1, 1}, j_{1, 2}, ..., j_{1, d}; j_{2, 1}, ....; j_{b, 1}, j_{b, 2}, ..., j_{b, d},$$
such that 
$$j_{k,\ell}\notin\Sigma_{k, \ell},$$
then  $(j_1, j_2, ..., j_b)\in\Cal D$.

We now show that the set defined in ($\diamondsuit\diamondsuit\diamondsuit$) is an infinite set for all $i$, by showing that 
it contains an infinite subset. Setting $i=1$, from Lemma 0,
the set 
$$F_{1, 1}:=\{x_1; x_1^2+1\in S_{qf}\}$$ is an infinite set. Since the non-admissible set in $x_1$,
$\Sigma_{1, 1}$ is finite, in fact $$|\Sigma_{1, 1}|\leq  (d+2)N(p, q, d),$$
the set in $x_1$
$$A_{1, 1}:=F_{1,1}\backslash \Sigma_{1, 1}$$
is infinite. Fix 
$$x_1\in  A_{1, 1},$$
and call $A_{1, 1}$ the admissible set (in $x_1$).
We may repeat the argument for $x_2$ by
defining 
$$F_{1, 2}:=\{x_2; x_2^2+x_1^2+1\in S_{qf}, \text{ fixed }x_1\in  A_{1, 1} \}.$$
Since $x_1^2+1$ is square-free by construction, Lemma 0 says that $F_{1, 2}$ is an infinite set.
Define the admissible set in $x_2$,
$$A_{1, 2}(x_1):=F_{1, 2}(x_1)\backslash \Sigma_{1, 2}(x_1),$$
which is again an infinite set. Fix 
$$x_2\in  A_{1, 2}.$$
$$\vdots$$
$$A_{1, d}(x_1, x_2, ..., x_{d-1}):=F_{1, d}(x_1, x_2, ..., x_{d-1})\backslash \Sigma_{1, d}(x_1, x_2, ..., x_{d-1})$$
is an infinite set and we fix 
$$x_d\in  A_{1, d}.$$
So $(\diamondsuit\diamondsuit\diamondsuit)$ holds for $i=1$. Clearly this construction maybe repeated for  
$i= 2, ..., b$ and concludes the proof. 
\hfill$\square$
\enddemo
\noindent{\it Remark 1.} The second part of the non-degeneracy condition (i) will only be used below to deal with the 
exceptional case when condition (iii) is not applicable. 

\noindent{\it Remark 2.} In fact, it follows as a Corollary of Theorem 1.1 in \cite{LX} that the good set has positive density in $(\Bbb Z^d)^b$.
We have kept the more elementary proof of a weaker assertion for the purposes here, since \cite{LX} is a rather involved paper
in analytic number theory.

\subheading {2.2 Size of connected sets on the characteristics} 
A set 
$$S\subseteq \Bbb Z^{b+d}$$
is called {\it connected}, if for all $a$, $b\in S$, there exist $a_1$, $a_2$, ..., $a_m\in S$, such that 
$$a_{k+1}-a_k\in\Gamma,\tag 2.10$$ for all $k\in\{0, 1, ... , m\}$ with $a_0:=a$, $a_{m+1}:=b$, 
where $\Gamma$ is as defined in (2.6). The number of elements in $S$, $|S|$, 
is its size. 

It follows that if $S$ is connected, then 
$$a_k-a_{k'}\in \tilde\Gamma_R,$$
for all $a_k$, $a_{k'}\in S$, $k\neq k'$, where $\tilde \Gamma_R$ as defined in ($\sharp$) for some $R>0$.
Choose an (arbitrary) element $a_0\in S$ and call it the root. Then 
\item{(P1)} $a_k-a_0\in \tilde\Gamma_R$, for all $k\neq 0$ 

and there must exist $k'\neq 0$, such that 
\item{(P2)} $a_{k'}-a_0\in \Gamma, \, a_{k'}\in S$.

We note that a subset of a connected set is {\it not} necessarily connected. 

Let $$\Cal S=\text{supp } u^{(0)}\subset \Cal C. \tag 2.11 $$
We consider the connected sets on (contained in) the characteristic  $\Cal C$. 
Below is the main result of the section:

\proclaim{Proposition 2.1}
Assume that $u^{(0)}=\sum_{k=1}^b a_k \cos ({-(\sqrt{j_k^2+1})t}+j_k\cdot x)$ is good satisfying the non-degeneracy conditions (i-iii). On the characteristic hyperboloid $\Cal C$, 
the connected sets are of size at most $\max (2d, 2b)$. If $b\geq d+1$, then 
the set $\Cal S$ is the only connected set of size $2b$, all the other connected sets  
are of size at most $2d$.    
\endproclaim

\demo{Proof}
If $(n, j)\in\Cal C$, then 
$$(n\cdot\omega^{(0)})^2-j^2-1=0.\tag 2.12$$ 
Therefore, as noted in (1.11),  
$n$ must be of the form $n=n_ke_k$ for some $k\in \{1, 2, ..., b\}$.
So the characteristic $\Cal C$ only consists of  ``singletons",
i.e., 
$$\Cal C\subseteq\{(n,j)|n=n_ke_k \text{ for some } k=1, 2, ..., b, n_k\in\Bbb Z\}.\tag \ddag$$
Combined with the definition of a connected set and its ensuing properties (P1, 2), it then
follows that (after designating a root) only $(\nu, \eta)$ of the form
considered in condition (iii) could possibly lead to connected sets on $\Cal C$. 

Assume that there is a connected set $S$ on the hyperboloid $\Cal C$. There are 
two cases: (a) For all $(n, j)\in S$, $|n|>p$, (b) there exists $(n, j)\in S$ with $|n|\leq p$. 

Case (a):  The size of $S$ must satisfy
$$|S|\leq 2.\tag 2.13$$
This is because if $|S|\geq 3$, then there must be a connected subset $S'\subseteq S\subset\Cal C$ satisfying $|S'|=3$. Let $(n, j)$, $(n', j')$ and $(n'', j'')$ be the 3 distinct points in $S'$, then 
$n$, $n'$ and $n''$ must be of the form $n=n_ke_k$, $n'=n'_ke_k$ and $n''=n''_ke_k$ 
for some $k=1, 2, ..., b$. Here 
we used that for $(\nu, \eta)\in \Gamma$, $|\nu|\leq p$ and property ($\dag\dag$). Call $(n, j)$ the root.

Let $$(\nu_1, \eta_1)=(n'-n, j'-j)=((n_k'-n_k)e_k, -(n_k'-n_k)j_k)$$ and 
$$(\nu_2, \eta_2)=(n''-n, j''-j)=((n_k''-n_k)e_k, -(n_k''-n_k)j_k).$$
Subtracting (2.12) evaluated at $(n', j')$ from the equation at $(n, j)$ and likewise (2.12) evaluated at $(n'', j'')$ from that at $(n, j)$ lead to a system of two linear equations in $(n, j)$:
$$\cases -2(n\cdot\omega^{(0)})(\nu_1\cdot\omega^{(0)})+2j\cdot\eta_1+\eta_1^2-(\nu_1\cdot\omega^{(0)})^2=0,\\
-2(n\cdot\omega^{(0)})(\nu_2\cdot\omega^{(0)})+2j\cdot\eta_2+\eta_2^2-(\nu_2\cdot\omega^{(0)})^2=0.\endcases\tag 2.14$$
After a straight forward computation,  they take the form (*) with 
$m=n_k$, $\ell=k$, $(m_\ell)_1=(m_k)_1=-n_k'$ and  $(m_\ell)_2=(m_k)_2=-n_k''$.
These two equations are incompatible if $n_k'\neq n_k''$, as observed previously in (2.7). So $|S|<3$. 

Case (b): Assume $|S|>2d$, then it must contain a {\it connected} subset $S'$, 
$$S'\subseteq S, \, |S'|=2d+1.$$  
Subtracting (2.12) evaluated at $(n', j')$ from the 
equation evaluated at $(n, j)$ for all $(n', j')\in S'\backslash \{(n, j)\}$, since $\eta=(n'-n, j'-j)\neq 0$ from condition (ii), there are $|S'|-1$ {\it proper} $d$-dimensional hyper-planes 
of the form in (*). If the non-degeneracy condition (ii) is applicable, then $|S'|\leq 2d$, which is a contradiction. 

If all $m_\ell=\pm m$ in the set $\sigma$ in condition (iii), we may assume $b\geq d+1$ as otherwise the intersection of $2d$ planes as in (iii) 
is empty by using (2.7, 2.8) and the conclusion of (iii) holds. In this case, the only possible solutions $(n, j)$ form the set 
$$M=\{(-me_k, mj_k); |m|=1, 2, ..., p \}_{k=1}^b, $$ using the second part of condition (i) and ($\ddag$). 
Since 
$$M\cap\Cal C=\Cal S,$$
the subset restricted to $|m|=1$, 
the only connected set of size $2b$ is the exceptional set $\Cal S$ defined in (2.11).  
\hfill$\square$
\enddemo

\noindent{\it Remark.} It is important to note that the system of linear equations in (2.14) are in the variables $n$ {\it and} $j$. 
After eliminating the variables $n$, it generally leads to {\it quadratic} polynomials in $j$, cf. sect.~5.1. (Here 
due to the very special property ($\ddag$), exceptionally, the system is linear in $j$.)
This is an essential complication compared to NLS, which mostly leads to linear systems in $j$, cf. sect.~2.3 of \cite{W2}. 

We have now achieved the block structure mentioned earlier, namely 
\proclaim {Corollary} If $u^{(0)}$ is good, satisfying the non-degeneracy conditions (i-iii), then  
restricting to $\tilde\Cal C:=\Cal C\backslash \text{supp } u^{(0)}$, where $\Cal C$ is the characteristic defined in (1.9) and $\text{supp } u^{(0)}$ in (1.7), the matrix $A$ defined in (2.5) is block-diagonal with blocks of size at most 
$2d\times 2d$, i.e.,  there is the block decomposition: $$A_{\tilde\Cal C}=\oplus_{\alpha} A_{\alpha},\tag$\bigstar\bigstar$ $$
where $\alpha$ are connected sets, and $A_\alpha$ are $A$ restricted to $\alpha$ -- therefore of size at most $2d\times 2d$.
\endproclaim 
\demo{Proof} This is an immediate consequence of Proposition~2.1 and the definition of {\it connected} in (2.10). \hfill$\square$
\enddemo

\head{\bf 3. Extraction of parameters -- solving the $Q$-equations}\endhead
We continue to work in the real; but for the sake of combinatorics, it is convenient to adopt complex 
notations. Let $$v^{(0)}=\sum_{k=1}^b a_k e^{-i(\sqrt{j_k^2+1})t}e^{ij_k\cdot x}$$ and ${\bar v}^{(0)}$
its complex conjugate. Then 
$$u^{(0)} =\frac {v^{(0)}+{\bar v}^{(0)}}{2}.$$

The nonlinear matrix equations in (2.2) are solved using the Lyapunov-Schmidt decomposition, as mentioned in sect.~1.7. 
Writing (2.2) as $$F(u)=0, $$ the $Q$-equations are the restrictions to the set $\Cal S$ defined in (2.11):
$$F(u)|_\Cal S=0;$$ the remaining equations are the $P$-equations. These equations are solved iteratively
using a Newton scheme similar to the one in \cite{W2}.

It is natural to start with the $Q$-equations and solve for the frequencies:
$$\omega_k=\sqrt{ j_k^2+1+\frac{1}{a_k}[\Cal M(\frac{v+{\bar v}}{2})^{*p+1}\Cal M] (-e_k, j_k)},$$
where $k=1, 2, ..., b$. For the first iteration, setting $u=u^{(0)}$, and in view of ($\flat$), we obtain
$$\omega_k^{(1)}=\sqrt{j_k^2+1}+\frac{1}{2^{p+2}a_k\sqrt{j_k^2+1}}(v^{(0)}+{\bar v}^{(0)})^{*p+1}{(-e_k, j_k)}+\Cal O(\delta^{2p}),$$
where $k=1, 2, ..., b$. So the frequency modulation: 
$$\aligned \Delta \omega_k^{(1)}:=\omega_k^{(1)}-\omega_k^{(0)}&=\frac{1}{2^{p+2}a_k\sqrt{j_k^2+1}}(v^{(0)}+{\bar v}^{(0)})^{*p+1}{(-e_k, j_k)+\Cal O(\delta^{2p}})\\
:&=\Omega_k+\Cal O(\delta^{2p}),\endaligned\tag 3.1$$
where $k=1, 2, ..., b$. 

There are the following estimates on amplitude-frequency modulation.
\proclaim{Proposition 3.1}  Assume that $u^{(0)}=\sum_{k=1}^b a_k \cos({-(\sqrt{j_k^2+1})t}+j_k\cdot x)$ is a solution to the linear equation with $b$ frequencies and 
$a=\{a_k\}_{k=1}^b\in (-\delta, \delta)^b\backslash\{0\}=\Cal B(0, \delta)=\Cal B\subset\Bbb R^b\backslash\{0\}$. Assume that $p$ is even and $\epsilon'\in (0,1)$. 
There exists a subset $\Cal B'\subset\Cal B$ with 
$$\text{meas }\Cal B' <\epsilon'\delta^{b}/2,$$ 
and $\delta_0>0$ such that 
if $a\in\Cal B\backslash\Cal B'$, an open set, and $\delta\in (0, \delta_0)$, then
$$\align 
&\Vert \Delta \omega^{(1)}\Vert\asymp \delta^{p},\\
&\Vert \frac{\partial \omega^{(1)}}{\partial a}\Vert \asymp \delta^{p-1},\\
&\Vert(\frac{\partial \omega^{(1)}}{\partial a})^{-1}\Vert \lesssim \delta^{-p+1},\\
&\big| \det(\frac{\partial \omega^{(1)}}{\partial a})\big| \gtrsim \delta^{(p-1)b},\endalign$$
where the constants implied by $\asymp$, $\lesssim$ and $\gtrsim$ only depend on $p, b, d$ and $\epsilon'$.  
\endproclaim

\noindent{\it Remark.} Note that $\omega=\omega^{(0)}+\Cal O(\delta^p)$ and $\omega^{(0)}$ is {\it fixed} -- this is why the initial 
``perturbation theory" in the forthcoming Lemmas~4.1 and 5.1 are {\it singular} and requires a {\it non-perturbative} treatment. Afterwards it returns
to {\it regular} perturbation theory in sect.~5.6. 

\demo{Proof}
Let $M$ be the convolution matrix:
$$M={(v^{(0)}*{\bar v}^{(0)}})^{*p/2}*,\tag 3.2$$
and $$D_k=2^{p+2}\sqrt {j_k^2+1}, \, k=1, 2, ..., b.$$
The proof is modelled after that of Proposition 3.3 in \cite{W2}. 
Using (3.1) and expanding the $(p+1)$-fold convolution, we obtain that 
$$\aligned B_k:=&D_k\Omega_k\\
=&C_{p+1}^{p/2} (M_{kk}+\sum_{i\neq k}\frac{M_{ki}a_i}{a_k}),\, k=1, ..., b,\endaligned\tag 3.3$$
cf. the second expression after (3.36) in \cite{W2}. Here we used that 
$${(v^{(0)}*{\bar v}^{(0)}})^{*p/2}*v^{(0)}$$
is the only contributing term. (This is the same term that appears in the NLS in (3.36) \cite{W2} and we shall use some of 
its properties derived there in (3.36)-(3.37).)

From the structure of $M$, 
$$M_{kk}(a_1, a_2, ..., a_b)=P(a_1, a_2, ..., a_b)$$
and $$M_{ki}=P'(a_1, a_2, ..., a_b)a_k\bar a_i,\, k\neq i$$
where $P$ and $P'$ are homogeneous polynomials in $a$ with positive integer coefficients
and are invariant under any permutations of the arguments, $P$ is of degree $p$, $P'$, $(p-2)$.
So $B_k$ is a homogeneous polynomial in $\{a_i\}_{i=1}^b$ of degree $p$ and can be written as 
$$B_k(a_1, a_2, ..., a_b)=C_{p+1}^{p/2} [P(a_1, a_2, ..., a_b)+P'(a_1, a_2, ..., a_b)P_k(\{a_i\})_{i\neq k}],$$
with $P_k=\sum_{i\neq k} a_i^2$ and $P$, $P'$ as above, for 
$k=1, 2, ..., b$.  

We first prove the last two estimates. Setting $a=(1,1, .., 1)$ and using that $P$, $P'$ and $P_k$ are polynomials 
with positive coefficients, we have
$$\frac{\partial B_k} {\partial a_i}(1, 1, ..., 1)> \frac{\partial B_k}{\partial a_k}(1, 1, ..., 1)$$
for all $i\neq k$. Let $q$ be the diagonal elements and $Q$ the off-diagonal ones at $(1, 1,..., 1)$. This
gives $q$, $Q\in\Bbb N^+$ satisfying 
$$0<q<Q.$$
For example, in the cubic case, $p=2$, $P=\sum_{i=1}^b a_i^2$, $P'=1$ and $P_k=\sum_{i\neq k}a_i^2$ 
giving $q=6$ and $Q=12$. 

Using the same argument as in the proof of Proposition 3.3, sect. 3.2  of \cite{W2},  the partial derivative matrix: 
$$[[\frac{\partial B_k}{\partial a_i}]]$$ has a simple eigenvalue $\lambda_1=q+(b-1)Q\neq 0$
and a $(b-1)$- fold degenerate eigenvalue $\lambda_2=q-Q\neq 0$. So 
$$\det \big (\frac{\partial B_k}{\partial a_i}\big)(1, 1, ..., 1)\neq 0=\det \big ( \frac{\partial B_k}{\partial a_i}\big)(0, 0, ..., 0).$$
Hence $\det \big(\frac{\partial B_k}{\partial a_i}\big)$ is not a constant.

Let $D$ be the $b\times b$ diagonal matrix with diagonals 
$D_k$, $k=1, 2, ..., b$. Since $\det \big(\frac{\partial B_k}{\partial a_i}\big)$ is a homogeneous polynomial in $a$ of degree 
at most $(p-1)b$ and $D$ is bounded and invertible, this proves 
the last two estimates taking into account also the $\Cal O(\delta^{p}) $ perturbation. 
As a consequence, this also proves the first two estimates.
\hfill $\square$
\enddemo

\noindent{\it Remark.} When $p$ is odd, $\Omega_k(a)=0$, $k=1, 2, ..., b$, for all $a$. There is no frequency modulation 
at order $\Cal O(\delta^p)$. As mentioned earlier, this is why $p$ is taken to be even. Note also that Proposition 3.1 does 
not need $u^{(0)}$ to be good, as by definition $\omega^{(1)}$ is {\it independent} of the correction $\Delta u^{(1)}$.

\head{\bf 4. The first step -- solving the $P$-equations}\endhead
Let $F'$ be the operator linearized at $u^{(0)}$ and evaluated at  $\omega=\omega^{(1)}$.
 Let $F'_N (\omega^{(1)}, u^{(0)})$ be the restricted operator:   
$$\aligned F'_N(n, j; n',j')&=F' (n, j; n',j'),\quad\text{if } (n, j), (n', j')\in [-N, N]^{b+d}\backslash\Cal S,\\
&=0,\qquad\qquad\qquad\,\,\text{otherwise,}\endaligned$$
where $\Cal S$ is as defined in (2.11).
In other words, the operator $F'_N$ is $F'$ restricted to the intersection of  the domaine of the $P$-equations with $[-N, N]^{b+d}$.
We make the first corrections to the unperturbed solution $u^{(0)}$ by 
solving the $P$-equations, using a Newton scheme.  The key is the invertibility of $F'_N$. For the first iteration,  we take $$N=|\log\delta|^s, \tag 4.1$$ 
for some $s>1$ to be determined in sect.~5.1, after the proof of Lemma~5.5, in (Fiv).

Since $$|n\cdot\Delta\omega^{(1)}|\leq \Cal O(|\log\delta|^s)\delta^p\ll c|\log\delta|^{-sq}, \quad c>0,\, s>1,\, q>1,$$
for small $\delta$, from (4.1) and (1.13), the resonance structure remains the {\it same} as for $\omega=\omega^{(0)}$ and one may
use the block-diagonal structure ($\bigstar\bigstar$) in the Corollary in sect.~2, to invert {\it block by block}. We obtain 
\proclaim{Lemma 4.1}
Assume that $u^{(0)}=\sum_{k=1}^b a_k \cos ( {-(\sqrt{j_k^2+1})t}+j_k\cdot x)$ is good satisfying the non-degeneracy conditions (i-iii) and $p$ even. Let $b>C_pd$, where $C_p$ only depends on $p$, and $\epsilon \in(0, 1)$. 
There exists $\delta_0>0$ such that for all $\delta\in (0, \delta_0)$, there exists a subset $\Cal B_{\epsilon',\epsilon, \delta}$, $(-\delta, \delta)^b\backslash\{0\}:=\Cal B\supset \Cal B_{\epsilon', \epsilon, \delta}\supset \Cal B'$,  the set in 
Proposition 3.1, with
 $$\text {meas }\Cal B_{\epsilon', \epsilon, \delta}<\epsilon'\delta^{b}.$$ On $\Cal B\backslash\Cal B_{\epsilon', \epsilon, \delta}$, an open subset,
the operator $F'_N$ satisfies:
$$\Vert [F'_N]^{-1}\Vert\leq \delta^{-p-\epsilon}, \tag 4.2$$
and there exists $\beta\in (0, 1)$, depending only on $\text {supp }u^{(0)}, p, b, d$ and the $H$ in (1.1)  such that 
$$|[F'_N]^{-1}(x, y)|\leq \delta^{\beta|x-y|}=e^{-\beta|\log\delta||x-y|}\tag 4.3$$
for all $x$, $y\in [-N, N]^{b+d}\backslash\Cal S$ such that $|x-y|>1/\beta^2$.
 \endproclaim

The upper bound in (4.2) is {\it non-perturbative}, since 
$\delta^{-p-\epsilon} \gg \delta^{-p}$, the latter being the inverse of smallness
of the perturbation $A$ in (2.5) of the linearized operator $F'$ in (2.3), and relies
fundamentally on the block structure. 
This decomposition also leads to a {\it geometric description} of the resolvent series used in deriving (4.3), see (4.11-4.16) below,
which could be of independent interest.

Let $P$ be the projection on $\Bbb Z^{b+d}$ onto $\Cal C$ defined in (1.9).
The following plays a key role toward proving Lemma~4.1.
\proclaim{Lemma 4.2} Let $u^{(0)}=\sum_{k=1}^b a_k \cos ( {-(\sqrt{j_k^2+1})t}+j_k\cdot x)$ be a good linear solutuion.
The linearized operator $F'$ evaluated at $\omega^{(1)}$, $u^{(0)}$, and restricted to $\tilde\Cal C=\Cal C\backslash \Cal S\subset \Bbb Z^{b+d}$: $PF'P$ can be written as 
$$PF'P=\oplus_\alpha\Cal F_\alpha,$$
where $\Vert \Cal F_\alpha\Vert =\Cal O(\delta^p)$, $\alpha$ are connected sets on $\tilde\Cal C$ satisfying 
$|\alpha|\leq 2d$. Moreover, 
$$\Cal F_\alpha= \text{diag }(\pm 2n\cdot\Omega\sqrt {j^2+1}\,)\big|_\alpha+A_\alpha,\tag 4.5$$
$\Omega=\{\Omega_k\}_{k=1}^b$ as defined in (3.1), $(n, j)\in\alpha$, $A_\alpha$ is  the $A$ defined in (2.5) restricted to $\alpha$, as in ($\bigstar\bigstar$).
Furthermore the (at most) $2d\times 2d$ matrix $\Cal F_\alpha$ satisfies $$\det \Cal F_\alpha\not \equiv  0,$$ if $b>C_{p, \alpha} d$. 
\endproclaim
\demo{Proof}  
Since $\omega^{(1)}$ only modifies the matrix entries and {\it not} the structure, the direct sum decomposition follows from the Corollary.
Equations (2.4, 3.1) are then used to arrive at the first matrix in (4.5).

To prove the second  part of the Lemma, we set
$$a_1=a_2= ...=a_b=1.$$ It is essentially a direct computation
similar to the proof of Lemma 4.2 in \cite{W2} and we shall use some 
of the derivations there.

We first note that if $(n, j)\in\Cal C$, then $n=n_ke_k$ for some $k=1, 2, ..., b$, and $n_k$, $j$ 
satisfy 
$$(n_k\omega^{(0)}_k)^2-j^2-1=0.$$
So 
$$\sqrt{j^2+1}=\pm n_k \sqrt{j_k^2+1}.$$
From (2.4) with $\omega=\omega^{(1)}$, (2.5, 3.1, 3.3) and using the above, it then follows that 
$$2^{p+1}\Cal F_\alpha(n, j; n, j)=-n_k^2 C_{p+1}^{p/2}(M_{11}+(b-1)M_{12})+(p+1)C_p^{p/2}M_{11},$$
where $M_{11}$, $M_{12}$ denote respectively the $(1,1)$, $(1,2)$ element of the matrix $M$  defined in (3.2)
and we used the symmetry: $M_{ii}=M_{11}$ for all $i$  and $M_{ij}=M_{12}$ for  all $i$, $j$, $i\neq j$.
$M_{11}$ is a polynomial in $b$ of degree $p/2$; while $M_{12}$, $p/2-1$. 

From  \cite{W2},  with $p\to p/2$, the combinatorial factor, the coefficient, in front of the $\Cal O(b^{p/2})$ term in $M_{11}$ is 
$$1+2!C_{p/2}^2+3!C_{p/2}^3+...+m!C_{p/2}^m+ ... + (p/2)!C_{p/2}^{p/2}; $$
while the coefficient in front of the leading order $\Cal O(b^{p/2-1})$ term for $M_{12}$ is 
$$C_{p/2}^1+2!C_{p/2}^2+3!C_{p/2}^3+...+m!C_{p/2}^m+ ... + (p/2)!C_{p/2}^{p/2}.$$
These two formulae appear just above (4.14) in \cite{W2} and are derived using the binomial expansion.

Write 
$$2^{p+1}\Cal F_\alpha(n, j; n, j)=Rb^{p/2}+\Cal O(b^{{p/2}-1}).$$
Using the above two formulae, it follows from direct computation that 
when $p=2$, $R=0$ if and only if $n_k=\pm 1$. When $p>2$, write
$N=n_k^2$. Setting $R=0$ leads to 
$$\aligned N&=\frac{(p/2+1)(\sum_{m=2}^{p/2} m!C_{p/2}^m+1)}{(p/2+1)+2\sum_{m=2}^{p/2} m!C_{p/2}^m}\\
&=\big(\frac{p+2}{4}\big)\big(\frac{1+\frac{1}{A}}{1+\frac{p+2}{4A}}\big),\\
&=\big(\frac{p+2}{4}\big)+\big(\frac{p+2}{4}\big)\big[\big(1+\frac{1}{A}\big)\sum_{n=1}^\infty(-x)^n+\frac{1}{A}\big]\endaligned$$
where $p>2$, $A=\sum_{m=2}^{p/2} m!C_{p/2}^m$ and $x=\frac{p+2}{4A}$. 

Since 
$$0<\big|\big(\frac{p+2}{4}\big)\big[\big(1+\frac{1}{A}\big)\sum_{n=1}^\infty(-x)^n+\frac{1}{A}\big]\big|<1/2 $$
from direct computation, 
$N\notin\Bbb Z$ for $p> 2$.
So $R\neq 0$ for integer $n_k^2$.  
Taking into account that $(\mp e_k, \pm j_k)$, $k=1, 2, ..., b$, are in $\Cal S$, this proves that  
$$\Cal F_\alpha (n, j; n, j)=\Cal O(b^{p/2}),$$ 
for all $(n, j)\in \alpha$ on $\Cal C\backslash \Cal S$.
Since all off-diagonal elements are of order at most $\Cal O(b^{p/2-1})$
and $\Cal F_\alpha$ is at most of size $2d\times 2d$, 
$$\det \Cal F_\alpha(1, 1, ..., 1)\neq 0$$ for $b>C_{p, \alpha}d$; so 
$$\det \Cal F_\alpha(1, 1, ..., 1)\not\equiv  0.$$
\hfill$\square$
\enddemo

\noindent{\it Remark.} In the proof below, Lemma 4.2 will be used only for a finite number ({\it independent} of $\delta$) of
blocks near the origin; for blocks away from the origin, one varies $\omega$. 

\subheading {4.1. Proof of Lemma 4.1} 
The proof is rather lengthy. It is therefore separated into two parts.
We first prove the norm estimate.
\demo{$\bullet$ Proof of (4.2) of Lemma 4.1}
One first noes that if $$\pm n\cdot\omega^{(0)}+\sqrt{j^2+1}\neq 0,$$ 
then $$|\pm n\cdot\omega^{(0)}+\sqrt{j^2+1}|\geq c'\Vert n\Vert_1^{-q},\tag 4.6$$ 
for some $c'>0$ and $q>b^2$, using \cite {Schm}, cf. also \cite{R} for the scalar case. 
This follows from multiplying the two expressions corresponding to the $\pm$ signs
in the absolute value and $|j|\lesssim |n|$ (otherwise the inequality in (4.6)
is trivially true). The bound in (4.6) implies that 
$$\Vert [P^cF'_NP^c]^{-1}\Vert\leq C'N^{q},\tag 4.7$$
where $P^c=I-P$,  for some $C'>0$ and small $\delta$.

From Schur's complement reduction \cite{S1, 2}, $\lambda$ is in the spectrum of $F'_N$, $\lambda\in\sigma(F'_N)$, if and only if
$0\in\sigma(\Cal H)$, where 
$$\Cal H=PF'_NP-\lambda-PF'_NP^c(P^cF'_NP^c-\lambda)^{-1}P^cF'_NP.\tag 4.8$$
Moreover (4.7) implies that (4.8) is analytic in $\lambda$ in the interval 
$$(-1/(2C'N^q), 1/(2C'N^q))$$ and on the same interval
$$\Vert PF'_NP^c(P^cF'_NP^c-\lambda)^{-1}P^cF'_NP\Vert\leq \Cal O(\delta^{2p}N^q).$$
So for small $\delta$, we only needs to prove invertibility of the first term in (4.8), which is a much smaller matrix.

Toward that purpose, we identify the set of connected sets 
$\{\alpha\}$ on $$\Cal C\cap [-N, N]^{b+d}\backslash \Cal S$$ 
with the set $\{1, 2,..., K_1\}$, where $K_1=K_1(N)$. So
$$PF'_{N}P=\oplus_k\Cal F_k(a),\quad k\leq K_1(N),$$
where each $\Cal F_k$ is of the form in (4.5).

Fix $$N_0=N_0 (p, b, d, \epsilon')\tag 4.9$$
large to be determined by (4.10) below.   
For a given $\Cal F_k$, define the support of $\Cal F_k$ to be 
$$\Bbb Z^{b+d}\times \Bbb Z^{b+d}\supset\text{supp }\Cal F_k=\{(x, y) |\Cal F_k(x,y)\neq 0\}.$$

For matrices $\Cal F_k$, such that $$\text{supp }\Cal F_k\cap \{[-N_0, N_0]^{b+d}\times [-N_0, N_0]^{b+d}\}\neq\emptyset,$$
we use the determinant. There are at most $K_0$ (independent of $\delta$) of these matrices.
Let  $$w=a\delta^{-1},$$ and 
$$P_k=P_k(a):=\det \Cal F_k(a\delta^{-1})=\det \Cal F_k(w)\not\equiv 0,$$
using Lemma 4.1. Therefore there exist $C_1$, $c_1>0$, such that given $\epsilon\in (0, 1)$, there exists $\delta_0\in (0, 1)$,  such that for all $\delta\in (0, \delta_0)$,
$$\text{ meas }\{a\in \Cal B||P_k|<\delta^{\epsilon}, \text{ all }k\leq K_0\}\leq C_1\delta^{b+c_1\epsilon}.$$
So $\Vert \Cal F _k^{-1}(a)\Vert\leq\Cal O(\delta^{-p-\epsilon})$ for all $k\leq K_0$. 

For matrices $\Cal F_k$ with $k>K_0$, 
$$\text{supp }\Cal F_k\cap \{[-N_0, N_0]^{b+d}\times [-N_0, N_0]^{b+d}\}=\emptyset$$
by definition. We use perturbation theory. This corresponds to case a) in the proof of Proposition~2.1.
So $\Cal F_k$ is at most a $2\times 2$ matrix and from ($\ddag$) there is $\ell\in\{1, 2, ..., b\}$ such that
$n\cdot \Omega=n_\ell\Omega_\ell$. Moreover from Lemma~4.2,
$$|\frac{\partial^2}{\partial\Omega_\ell^2}P_k|+|\frac{\partial}{\partial\Omega_\ell}P_k|>\frac{1}{2}\tag 4.10$$
for $N_0$ large enough depending only on $p$, $b$, $d$ and $\epsilon'$, where we used the form of the matrix in (4.5)  and
$$\Vert \frac{\partial \Omega}{\partial w}\Vert\asymp \Vert\big( \frac{\partial \Omega}{\partial w}\big)^{-1}\Vert\asymp \Cal O(1),$$
from Proposition~3.1 with the constants implied by $\asymp$ depending only on $p, b, d$ and $\epsilon'$, and that there are only finite types of ``convolution" matrices $A_k$.

Since $\Vert \Cal F_k(w)\Vert \leq \Cal O(|\log\delta|^{2s})$, (4.10) gives $\Vert \Cal F_k(w)^{-1}\Vert~\leq ~\Cal O(\delta^{-\epsilon})$
for all $K_0<k\leq K_1(N)$ away from a set in $w$ of measure less than $\epsilon'/2+\delta^{\epsilon/3}$, where we also
used $K_1(N)\leq\Cal O(|\log\delta|^{2(b+d)s})$, $s>1$.  
So $\Vert \Cal F_k(a)^{-1}\Vert~\leq ~\Cal O(\delta^{-p-\epsilon})$
away from a set in $a$ of measure less than $(\epsilon'/2+\delta^{\epsilon/3})\delta^b$.

Combining the above two regions,  one has that away from a set in $a$ of measure less than $\epsilon'\delta^{b}$,
$$\Vert [PF'_N(a)P]^{-1}\Vert \leq  \Cal O(\delta^{-p-\epsilon}).$$
The Schur reduction in (4.8) then gives (4.2), cf. Lemma 4.8 in \cite{BGS}. 
\hfill$\square$
\enddemo
\demo{$\bullet$ Proof of (4.3) of Lemma 4.1} Let $w=a\delta^{-1}$ as before and
$$\aligned \tilde F&=\oplus_\alpha \delta^{p}\Cal F_{\alpha}\oplus \text{ diag }[-(n\cdot {\omega^{(1)}})^2+j^2+1+\delta^{p}A(n,j;n,j)] |_{(n,j)\notin\Cal C}\\
:&=\oplus_\alpha \delta^{p}\Cal F_{\alpha}\oplus \Cal D_{\Bbb Z^{b+d}\backslash\Cal C}\endaligned\tag 4.11$$
where the first direct sum is exactly as in Lemma 4.2, with $\alpha$ connected subsets of $\Cal C$, $\tilde F:=\tilde F(w)$, $\Cal F_\alpha:=\Cal F_\alpha(w)$ and $A:=A(w)$ as in (2.5). 

To obtain the point-wise estimates, we use a resolvent expansion about $\tilde F$. 
For simplicity, write $x$ for $(n, j)$ etc. The matrix $A$ is a convolution matrix with diagonal $\bar A$. Let $\Lambda$ be the matrix: 
$$\aligned \Lambda(x, y)=\Lambda(y, x)&=(A-\bar A I)(x, y),  \text{ if } (x, y)\notin \Cal C\times\Cal C \\
&=0\qquad\qquad\qquad\, \text{ otherwise},\endaligned\tag 4.12$$
$\Vert \Lambda\Vert =\Cal O(\delta^p)$. 
We have
$$F'=\tilde F+\Lambda=\tilde F+\Lambda.$$
The resolvent expansion gives:
$$[F_N']^{-1}=[\tilde F_N]^{-1}-[\tilde F_N]^{-1}\Lambda_N[\tilde F_N]^{-1}+[\tilde F_N]^{-1}\Lambda_N[\tilde F_N]^{-1}\Lambda_N[F_N']^{-1},\tag 4.13$$
where as before the sub-index $N$ denotes the restriction to $[-N, N]^{b+d}$.  The analysis of the series is similar to the proof of (3.3) in 
Lemma 3.1 in \cite{W2}. Below we summarize the key steps.
For simplicity of notation, the subscript $N$ is omitted. 

We estimate the second term in the expansion:  
$$[\tilde F]^{-1}\Lambda[\tilde F]^{-1}.\tag 4.14$$


We need to estimate
$[\Cal F_\alpha]^{-1}\Lambda[\Cal F_{\alpha'}]^{-1}$,  $[\Cal F_\alpha]^{-1}\Lambda[\Cal D]^{-1}$, $[\Cal D]^{-1}\Lambda[\Cal F_{\alpha'}]^{-1}$ and $[\Cal D]^{-1}\Lambda[\Cal D]^{-1}$.
From the definition of $\Lambda_1$ in (4.12), the first term
$$[\Cal F_\alpha]^{-1}\Lambda[\Cal F_{\alpha'}]^{-1}=0.$$
Using (4.2) and since 
$$\Vert\Cal D^{-1}\Vert \leq \Cal O(|\log\delta|^{sq})$$
from (4.6) and small $\delta$,  summing over the last three terms yields 
$$\Vert [\tilde F]^{-1}\Lambda[\tilde F]^{-1}\Vert \leq \Cal O(\delta^{p-\tilde\epsilon}) \Vert [\tilde F]^{-1}\Vert,\,  0<\tilde\epsilon<1.\tag 4.15$$  

Iterating the resolvent expansion in (4.13) $r$ times yields the $(r+1)$ term series 
$$[F']^{-1}=[\tilde F]^{-1}-[\tilde F]^{-1}\Lambda[\tilde F]^{-1}+...+(-1)^r[{\tilde F}^{-1}\Lambda]^r[F']^{-1}.\tag 4.16$$
We note that the blocks in $\tilde F$ (and hence $[\tilde F]^{-1}$) are of sizes at most $2d$ and that 
$\Lambda$ satisfies $\Lambda(x, y)=0$ if $|x-y|>C$, for some $C$ depending only on 
 $\text{supp } u^{(0)}$, $p$, $b$ and $d$  in (1.1). For all $x, y$, matrix multiplication then infers that
in the resolvent series (4.16), the first 
$$r=\frac{|x-y|}{2Cd}$$ 
(after taking the integer part) terms are identically zero.

There exists $\beta\in (0, p/4Cd)$, such that for all $x, y$, such that $|x-y|>1/\beta^2$,  
iterating the bounds in (4.15) and using (4.2) to estimate the last, the $(r+1)$th term 
produces (4.3).
\hfill$\square$
\enddemo
\subheading{4.2. The first approximate solution} 
The nonlinear matrix equation (2.2) is now ready to be solved iteratively. The solution $u$ is held fixed on 
$\Cal S$: $$u(\mp e_k, \pm j_k)=a_k/2,\, k=1, ... , b.$$
The set $\Cal S$ is the domain of the $Q$-equations and is used to solve for the frequencies. The $Q$-equations were solved in sect.~ 3. 
To solve the $P$-equations, restrict the domain to 
$$[-N, N]^{b+d}\backslash \Cal S, \text{ where } N=|\log\delta|^s, s>1,$$
for the first iteration
and define $$\Delta u^{(1)}:= u^{(1)}-u^{(0)}=-[F'_N]^{-1}(\omega^{(1)}, u^{(0)})F(\omega^{(1)}, u^{(0)}).$$
We precipitate that the {\it quadratic} Newton scheme shall {\it compensate} the ``bad" $[F'_N]^{-1}$ estimate in (4.2), by giving the ``good" 
estimate on $F$ in (4.18) below, thus weld together singular perturbation theory with the regular perturbation theory starting in sect.~5.3.

Below is a summary of the precise findings.
\proclaim{Proposition 4.3}  Assume that $u^{(0)}=\sum_{k=1}^b a_k \cos (-(\sqrt {j_k^{2}+1})t+j_k\cdot x)$ a solution to the linear equation with $b$ frequencies is good and 
$a=\{a_k\}\in (-\delta, \delta)^b\backslash\{0\}=\Cal B\subset\Bbb R^b\backslash\{0\}$. Let $b>C_pd$, where $C_p$ only depends on $p$, and $\epsilon', \epsilon\in (0,1)$.
There exists $\delta_0>0$, such that
for all $\delta\in (0,\delta_0)$, there is a subset $\Cal B_{\epsilon', \epsilon, \delta}\subset\Cal B$ with 
$$\text{meas }\Cal B_{\epsilon', \epsilon, \delta} <\epsilon'\delta^{b}.$$
Let $\rho$ be a weight on $\Bbb Z^{b+d}$ satisfying 
$$\aligned\rho(x)=&e^{\beta|\log\delta| |x|},\, 0<\beta<1\, \text{ for } |x|>1/\beta^2,\\ 
=&1,\qquad\qquad\quad\,\qquad \quad\,\text{  for } |x|\leq 1/\beta^2.\endaligned$$
Define the weighted $\ell^2$ norm:
$$\Vert\cdot \Vert_{\ell^2(\rho)}=\Vert\rho\cdot\Vert_{\ell^2}.$$

There exists $\beta\in(0, 1)$, determined only by  
$\text{supp } u^{(0)}$, $p$, $b$, $d$ and $H$ in (1.1), such that if $a\in\Cal B\backslash\Cal B_{\epsilon', \epsilon, \delta}$, an open subset,
then 
$$\align&\Vert \Delta u^{(1)}\Vert _{\ell^2(\rho)}\lesssim \delta^p,\tag 4.17\\
&\Vert F(u^{(1)})\Vert_{\ell^2(\rho)}\lesssim \delta^{3p-1},\tag 4.18\\
&\Vert \Delta \omega^{(1)}\Vert\asymp \delta^{p},\tag 4.19\\
&\Vert \frac{\partial \omega^{(1)}}{\partial a}\Vert\asymp \delta^{p-1},\tag 4.20\\
&\Vert(\frac{\partial \omega^{(1)}}{\partial a})^{-1}\Vert\lesssim \delta^{-p+1},\tag 4.21\\
&\big| \det(\frac{\partial \omega^{(1)}}{\partial a})\big|\gtrsim \delta^{(p-1)b},\tag 4.22\endalign$$
where $\omega^{(1)}$ and $\Delta \omega^{(1)}$ as defined in (3.1); the implied constants in (4.17, 4.21, 4.22) depend on $p, b, d$ and $\epsilon'$ ;
while that in (4.18-4.20), only on $p, b$ and $d$.
Moreover $\omega^{(1)}$ is Diophantine
$$\Vert n\cdot \omega^{(1)}\Vert_{\Bbb T}\geq\frac{\xi}{|n|^\gamma},\quad n\in[-N, N]^b\backslash\{0\},\,\xi>0, \gamma>2b,\tag 4.23$$
where $\Vert\,\Vert_{\Bbb T}$ denotes the distance to integers in $\Bbb R$, $\xi$ and $\gamma$ only depend on $\omega^{(0)}$.
\endproclaim
\demo{Proof} We only need to prove (4.17, 4.18, 4.23). The rest is the content of (proven) Proposition 3.1.
Since $$F(u^{(0)})=\hat \Cal M(u^{(0)})^{*(p+1)}\hat \Cal M=(u^{(0)})^{*(p+1)}.$$
The support of $F$, $\text{supp }F=\{(\nu, \eta)\}$, consists of 
$(\nu, \eta)$ of the form 
$$(\nu, \eta)=(-\sum_{k=1}^b m_ke_k, \sum_{k=1}^b m_kj_k),$$
with $\sum_{k=1}^b|m_k|\leq p+1$.
Using the non-degeneracy condition (ii) (property ($\dag\dag$)) and that $$(-\ell e_k, \ell j_k)\notin\Cal C, \, k=1, 2, ..., b,$$ for $|\ell|\neq 1$,  it follows that 
$$F(u^{(0)}) \cap \Cal C\backslash \Cal S=\emptyset.$$ 

From the Newton scheme
$$\aligned \Delta u^{(1)}&=-[F_N'(u^{(0)})]^{-1}F(u^{(0)})\\
&=-{\Cal D}^{-1}F+(F'_N)^{-1}(F'_N-\Cal D){\Cal D}^{-1}F\\
&\sim |\log\delta|^{sq}\delta^{p+1}(1+\delta^{-\epsilon})\\
&=\Cal O(\delta^p)\endaligned$$
in $\ell^2$ norm, where we used (4.2), $\Cal D$ as defined in (4.11), we used (4.6) and small $\delta$.
Using the point-wise estimate (4.3), the (at least) exponential off-diagonal decay of $F'_N$, and since $\text{supp }F$ is a compact set , the above bound remains valid in the weighted $\ell^2$ norm, 
$\Vert\cdot\Vert _{\ell^2(\rho)}$ and we obtain (4.17).

Since the $Q$-equations are solved exactly, $F_{\Cal S} (u+\Delta u)=0$. (For simplicity, the superscripts have been dropped.) To prove (4.18), 
we only need to be concerned with  $F_{\Cal S^c} (u+\Delta u)$. Below (and in general) for notational simplicity, we omit the subscript $\Cal S^c$ and write 
$$\aligned F(u+\Delta u)&= F(u)+F'(u)\Delta u + \Cal O (\Vert F''(u)\Vert \Vert \Delta u\Vert^2),\\ 
&=-(F'-F'_N)[F'_N]^{-1}F(u)+\Cal O (\Vert F''(u)\Vert \Vert \Delta u\Vert^2).\endaligned$$
Since $$\aligned [F'-F_N'](x, y)=&0,\qquad  \qquad\quad  x, y\in [-N, N]^{b+d}\backslash \Cal S,\\
=&F'(x, y),\qquad  \text{otherwise},\endaligned$$
using the exponential off-diagonal decay of $F'$, (4.3), compactness of $\text{supp }F$
and (4.17), we obtain
$$\Vert F(u+\Delta u)\Vert_{\ell^2(\rho)}=\Cal O(\delta^{3p-1}).$$

The (linear) Diophantine property in (4.23) is a consequence of ($\dag$) in sect.~2.1,  using the 
Diophantine approximation result in \cite{Schm}, 
small $\delta$ and (4.1).
\hfill $\square$ 
\enddemo
\smallskip
\noindent{\it Conventions on constants}

Since $d$, $p$, $b>C_p d$, $\text{supp }u^{(0)}$ and the $\Cal M$ in (1.1) are {\it fixed},
from now on, constants which only depend on them will be denoted generically by $\Cal O(1)$; more generally,
constants which depend on fixed parameters will be denoted by $\Cal O(1)$ below.

\head{\bf 5. The $\theta$ estimates}\endhead
Proposition~4.3 puts the construction in a non-resonant form with $\omega^{(1)}$ as the 
parameter.  It provides the input for the initial scales in the Newton scheme. To 
continue the iteration, we need the analogues of Lemma~4.1 and Proposition~4.3 at larger scales. 
This section prepares the way toward that goal by proving Lemma~5.1 and Proposition ~5.4. 

Let $T=F'$ be the linearized operator defined as in (2.3-2.5) and the restricted operator 
$T_N=F'_N$ as defined above (4.1). 
To increase the scale from $N$ to a larger scale $N_1$, we pave the $N_1$
cubes with $N$ cubes. We add a one dimensional parameter $\theta\in\Bbb R$  and consider $T(\theta)$: 
$$\aligned
T(\theta) =&\text{diag }[-(n\cdot\omega+\theta)^2+{j^2+1}]+A\\
=&D'(\theta)+A\\
\endaligned\tag 5.1$$
where $\omega$, standing for $\omega^{(1)}$, is considered as a parameter in this section, $A$ defined as in (2.5) 
As before, the matrix:
$$A=A(a, \omega, u)$$ is T\"oplitz in the $n$-direction and {\it independent} of $\theta$.  
We make estimates in $\theta$ in this section. 

The one dimensional parameter $\theta$ is added, because $n$ and $\omega$ only appear 
as $n\cdot\omega\in\Bbb R$ on the diagonal. It is an {\it auxiliary} variable to facilitate the analysis. 
Using the covariance of $n\cdot\omega+\theta$, all estimates in $\theta$ are transformed into estimates 
in $\omega$ in the Newton construction of $u$, and $\theta$ is always {\it fixed} at $0$ in sect.~6.

In sects.~5.1-5.4, we derive estimates on $T^{-1}(\theta)$ for initial scales. This step is {\it non-perturbative}
 (cf. (5.2, 5.4) below), and is the $\theta$-analogue of the results in sects.~2 and 4. 
 In sects.~5.6 and 5.7, we iterate to obtain estimates on 
$T^{-1}(\theta)$ for all scales. Sect.~5.5 gives a general proof of Diophantine $\omega$, which is needed 
for the iterations starting in sect.~5.6 (in fact for any KAM-type iterations). 

\subheading{5.1 The initial estimate in $\theta$}
  Let $N=|\log\delta|^s$ $(s>1)$ be as in Proposition 4.3, and $T_N(\theta)=T_N(\theta; u^{(1)})$.
We first state the estimates (the $\theta$-analogue of Lemma 4.1), which will be proved in the course 
of sects.~5.2-5.4.
\proclaim{Lemma 5.1}
Let $u^{(0)}=\sum_{k=1}^b a_k \cos ( {-(\sqrt{j_k^2+1})t}+j_k\cdot x)$ be 
a solution to the linear Klein-Gordon equation (1.2) satisfying the non-degeneracy conditions (i, ii), and (4.17, 4.19) hold
with $\epsilon\in (0, 1/2)$.
Let $\sigma, \tau$ be numerical constants satisfying 
$$0<\tau<1/s<\sigma<1.\tag $\spadesuit$ $$ 
There exists $\delta_0>0$, such that for all $\delta\in(0, \delta_0)$
$$\Vert [T_N(\theta)]^{-1}\Vert\leq \delta^{-p-\epsilon}<e^{N^\sigma}, \tag {5.2}$$
and there exists $\beta\in (0,1)$, 
such that
$$|[T_N(\theta)]^{-1}(x,y)|\leq \delta^{\beta|x-y|}=e^{-\beta|\log\delta||x-y|}\tag 5.3$$
for all $x, y$ such that $|x-y|>1/\beta^2$, for $\theta$ away from a set $B_N(\theta)\subset\Bbb R$ with 
$$\text{meas }B_N(\theta)<\delta^{p+\epsilon/8b}<e^{-N^\tau}.\tag 5.4$$
\endproclaim 
The proof of Lemma 5.1 is related to that of Lemma 4.1. A non-perturbative approach is necessary here, 
because even though the first variation in $\theta$ of each diagonal element in $D$ near its two zeroes is 
of order $\Cal O(1)$, perturbation about each diagonal to achieve invertibility leads to
excisize a set in $\theta$ of measure at least   
$$\Cal O(|\log\delta|^{s(b+d)})\delta^p\gg\Cal O(\delta^p),$$
where $\Cal O(|\log\delta|^{s(b+d)})$ is the volume factor. Since the
estimate in $\theta$ will be translated into estimates in $\omega$ and the variation
in $\omega$ is only of order $\delta^p$, the above is too large to be useful.

So the proof again rests on variable reductions to achieve 
a block structure, and to invert block by block, 
as in the proof of Lemma 4.1.
However, due
to the presence of $\theta$,  (1.11) is not applicable and the geometric non-degeneracy condition (iii) cannot be used
-- only (i) and (ii) are at our disposal.  
The advantage of this exclusive reliance on the arithmetic condition (ii) is that the new block structure
(alternatively the connected set structure) revealed by Lemma~5.2 below holds in greater generality. 
But the variable reduction also becomes more difficult 
than that in sect.~2. The proof is divided into three steps, realized in sects.~5.2-5.4. 

\subheading{5.2 Spacing of zeroes} 
On the diagonal of $T_N$, sits the family of quadratic polynomials in $\theta$:
$$P(\theta)=-(n\cdot\omega^{(0)}+\theta)^2+{j^2+1},\, (n, j) \in [-N, N]^{b+d}.$$
Similar to the proof of Lemma 4.1, the zeroes of the polynomials, i.e., $\theta$ such that
$$P(\theta)=-(n\cdot\omega^{(0)}+\theta)^2+{j^2+1}=0,\, (n, j) \in [-N, N]^{b+d},$$
play an essential role.  These quadratic polynomials yield the roots
$$\Theta=-n\cdot\omega^{(0)}\pm \sqrt{j^2+1}, \, (n, j) \in [-N, N]^{b+d}.\tag 5.5$$ 
We note that for each {\it given} $(n, j)$, the roots are {\it simple}. 

The new element in the proof, compared to that of Lemma 4.1, is 
the spacing of the $\Theta$'s, i.e., the spacing of square roots. 
(In Lemma 4.1, $\Theta$ is restricted to $\Theta=0$ only.) 
Fix $i=(n, j)$ and $i'=(n', j')$, $i\neq i'$.
Denote by $\Theta_i$ and $\Theta_{i'}$ ($i\neq i'$), two roots defined in (5.5). 
Define $\rho$ to be the spacing between them: 
$$\rho:=\Theta_i-\Theta_{i'}.$$
Clearly $\rho$ takes one of the four possible forms below: 
$$\aligned\rho:=&\Theta_i-\Theta_{i'}\\
=&(n-n')\cdot\omega^{(0)}\pm\sqrt{j^2+1}\pm\sqrt{{j'}^2+1}\\
=&\nu\cdot\omega^{(0)}\pm\sqrt{j^2+1}\pm\sqrt{{j'}^2+1},\endaligned$$
where $\nu=n-n'\in[-2N, 2N]^b$. We have the following dichotomy:

\proclaim{Lemma D} Let $N=|\log\delta|^s$ ($s>1$), $\rho$ and $\nu$ as defined above.
\item{(D1)} If $\rho=0$, then $\nu$ has at most $2$ non-zero components.
\item{(D2)} If $\rho\neq 0$, then $$|\rho|\geq \frac{1}{|\log\delta|^{s\Cal L}}\quad (\nu\neq 0), $$
for some $\Cal L>1$. 
\endproclaim

\noindent{\it Remark.} Here it is essential that the constant in the upper bound in (D2) is {\it independent} of $j, j'\in [-N, N]^d$;
for small $\delta$, it can be chosen to be $1$ by increasing $\Cal L$. 

\demo{Proof} One may always write 
$$\sqrt{j^2+1}=m_1\tau_1, \, \sqrt{{j'}^2+1}=m_2\tau_2,$$
where $m_1, m_2\in\Bbb N$, $\tau_1$ and $\tau_2$ are square roots of square-free integers or the integer $1$. 
Clearly if $\nu$ has at least $3$ non-zero components, then 
$$\rho=\nu\cdot\omega^{(0)}\pm m_1\tau_1\pm m_2\tau_2$$
is a linear form with at least $1$ square root of square-free integer.
So $$\rho=\nu\cdot\omega^{(0)}\pm m_1\tau_1\pm m_2\tau_2\neq 0$$
follows from algebraic independence, \`a la ($\dag$), proving (D1). 

Below for concreteness, assume that $\rho$ is of the form  
$$\rho=\nu\cdot\omega^{(0)}+\sqrt{j^2+1}-\sqrt{{j'}^2+1}.$$
(The other three forms work the same way.)
If $\rho\neq 0$, one may assume that 
$$\rho_1:=\nu\cdot\omega^{(0)}-(\sqrt{j^2+1}-\sqrt{{j'}^2+1})\neq 0;$$
$$\rho_2:=\nu\cdot\omega^{(0)}+(\sqrt{j^2+1}+\sqrt{{j'}^2+1})\neq 0;$$
$$\rho_3:=\nu\cdot\omega^{(0)}-(\sqrt{j^2+1}+\sqrt{{j'}^2+1})\neq 0.$$
As otherwise 
if $\rho_1=0$, then 
$$\nu\cdot\omega^{(0)}=\sqrt{j^2+1}-\sqrt{{j'}^2+1};$$
if $\rho_2=0$, then 
$$\nu\cdot\omega^{(0)}=-(\sqrt{j^2+1}+\sqrt{{j'}^2+1});$$
and if $\rho_3=0$, then 
$$\nu\cdot\omega^{(0)}=\sqrt{j^2+1}+\sqrt{{j'}^2+1}.$$
Substituting into the expression for $\rho$ gives
$$|\rho|\geq 2\min(|\nu\cdot\omega^{(0)}|, 1)\geq  \frac{1}{|\log\delta|^{s\Cal L}}\quad (\nu\neq 0),$$
using \cite{Schm}.
So (D2) is satisfied. 

Multiplying $\rho$ by $\rho_1$, $\rho_2$ by $\rho_3$ produce
$$(\nu\cdot\omega^{(0)})^2-(j^2+{j'}^2+2)\neq \pm 2\sqrt{j^2+1}\sqrt{{j'}^2+1}.$$
Squaring yields 
$$I:=(\nu\cdot\omega^{(0)})^4-2(\nu\cdot\omega^{(0)})^2(j^2+{j'}^2+2)+(j^2-{j'}^2)^2\neq 0.\tag \dag\dag\dag$$
Rewrite the difference of the first two terms in $I$ 
as 
$$\sum_{k'=1}^{b'} C_{k'}w_{k'}+N',$$
where $w_{k'}\neq 1$ are square roots of square-free integers determined by the prime decompositions 
of $(j_k^2+1)$, $k=1, 2, ..., b$, $C_{k'}$, $N'\in \Bbb Z$. If all $C_{k'}=0$, then since
$I\neq 0$, 
$$|I|\geq 1;$$
otherwise using \cite{Schm} and the restrictions on $\nu$, $j$, $j'$ yields (D2). Evidently the other three
possibilities for $\rho$ yields (D2) as well and concludes the proof.
\hfill$\square$
\enddemo

\noindent{\it Remark.} We note that unlike ($\dag\dag$), ($\dag\dag\dag$) is weaker than a Diophantine
property as it is only for certain rational combinations, but it suffices to establish Lemma 5.1.

\subheading{5.3 Size of connected sets on the characteristics}
Let $\Cal Z$ be the set of $\Theta$ defined in (5.5). For each $\Theta\in \Cal Z$, define
the characteristic to be 
$$\Cal C(\Theta)=\{(n, j)\in\Bbb Z^{b+d}|-(n\cdot\omega^{(0)}+\Theta)^2+j^2+1=0\};\tag $\natural$ $$
and the two branches
$$\Cal C_{\pm}(\Theta)=\{(n, j)\in\Bbb Z^{b+d}|\pm(n\cdot\omega^{(0)}+\Theta)+\sqrt{j^2+1}=0\}.$$
Define the connected sets on $\Cal C(\Theta)$ as in the paragraph containing (2.10), at the very beginning of sect. 2.2. 

\proclaim{Lemma 5.2} Let $u^{(0)}=\sum_{k=1}^b  a_k \cos (-(\sqrt{j_k^2+1} )t + j_k\cdot x)$
satisfy the non-degeneracy conditions (i, ii). Then 
the connected sets $\alpha$ on $\Cal C(\Theta)$ are at most of size $4b$ for all 
$ \Theta \in \Cal Z$.
\endproclaim
\demo{Proof} We work separately on $\Cal C_{+}$ and $\Cal C_-$. 
If there are two distinct points $(n, j)$, $(n', j')\in\Cal C_+$:
$$\cases(n\cdot\omega^{(0)}+\Theta)-\sqrt{j^2+1}=0,\\
(n'\cdot\omega^{(0)}+\Theta)-\sqrt{{j'}^2+1}=0,\endcases$$
and are, moreover,  connected, then 
$$(n'-n, j'-j)=(\nu, \eta)\in \Gamma$$ by definition,
where $\Gamma$ as defined in (2.6). 
Since $\nu\neq 0$, if $\eta=0$, these two equations are 
incompatible using ($\dag$). One may therefore assume below that $\eta\neq 0$.
Subtracting the second from the first and squaring lead to the following equation:
$$2j\cdot\eta+\eta^2-(\nu\cdot\omega^{(0)})^2-2(\nu\cdot\omega^{(0)})\sqrt{j^2+1}=0,\tag **$$
cf. (2.14). As before, call $(n,j)$ the root. 

More generally, if there is a connected set of $(r+1)$ sites on $\Cal C_+$, choosing a 
root, there is a system of $r$ equations of the form (**), with 
$(\nu, \eta)\in \tilde\Gamma$, where $\tilde\Gamma$ is as defined in ($\sharp$).   Let $S=\{(\nu_i, \eta_i)\}_{i=1}^r$ be the set 
of $(\nu, \eta)$ appearing in (**). 
From the dichotomy (D1), if 
$$
(\nu_i, \eta_i)\in S, \, i=1, 2, ..., r,$$
then $\eta_i$ is a function of at most $2$ variables in $\{j_k\}_{k=1}^b$
and so are the difference functions 
$$\eta_{i}-\eta_{i'}$$
for all $i$, $i'=1, 2, ..., r$, $i\neq i'$. 
Moreover one may assume that $\eta_i\neq 0$ and $\eta_i-\eta_{i'}\neq 0$ 
for all $i$, $i'$, $i\neq i'$, as mentioned earlier. (Otherwise there is a contradiction to $S$ being connected.) 

Up to permutations of the set $\{1, 2, ..., b\}$, it follows that $S$ takes two possible forms
(for simplicity, we omit the sub-index on $\nu$ and $\eta$):
\item {(i)} $$S=\{(\nu,\eta)\},\, |S|=3,$$
$$\{\eta\}=\{\alpha_1j_1+\alpha_2j_2, \alpha_2j_2+\alpha_3j_3, \alpha_3j_3+\alpha_1j_1; \alpha_i\neq 0, i=1,2, 3\}.$$ 

\item {(ii)} $$S=\{(\nu,\eta)\},\, |S| \text{ a priori arbitrary},$$
and $\eta$ of form:
$$\aligned &\alpha_1'j_1,\\
&\alpha_1j_1+\alpha_2j_2,\\
&\alpha_1j_1+\alpha_3j_3,\\
&\quad \cdots\\
&\alpha_1j_1+\alpha_bj_b,\endaligned$$
where $\alpha_1'$, $\alpha_i\in\Bbb Z$, $i=1, 2, ..., b$. 

We proceed to bound the $S$ in (ii) by contradiction. The idea is to take advantage of the restriction on the variable dependence of $\eta$
to reach a system of equations of the form (**), which is ``not too large" and {\it must} contain a ``degenerate" subsystem, i.e., there exist $\alpha_i\neq 0$, for all $i$,
such that $$\sum_i \alpha_i\eta_i\equiv 0,$$
for the subsystem. Using ($\bigstar$), this leads to 
$$\sum_i \alpha_i\nu_i\equiv 0.$$
So for the degenerate subsystem, this specific linear combination of the third (nonlinear) term in (**) disappears! As it turns out, to deal with the degenerate subsystem, we only need to concern with the $r=2$ and $3$ 
systems. So below we start with them.
  
Assume that $r=2$ and that 
 there is a constant $C\neq 0$ such that $\eta_2 \equiv C\eta_1$, $\nu_2=C\nu_1$. From (D1) and without loss of generality,
one may assume that 
$$0\neq\eta_1=mj_1+hj_2.$$
Then 
$$\eta_1^2-(\nu_1\cdot\omega^{(0)})^2=-m^2-h^2+2mh(j_1\cdot j_2-\sqrt{j_1^2+1}\sqrt{j_2^2+1})\neq 0,$$
for $j_1$, $j_2\in\Bbb Z^d\backslash\{0\}$. 
So the two equations of the form (**) are incompatible if $C\neq 1$. Therefore $|\alpha|\leq 2$. (Geometrically, this is because of curvature, as for (2.7, 2.8).)

Assume that $r=3$ and that 
there exist constants $C_1\neq 0$, $C_2\neq 0$ such that 
$$C_1\eta_1+C_2\eta_2+\eta_3\equiv 0,$$
and there are {\it no} constants $C'_i$, such that 
$$C'_i\eta_i\equiv \eta_{i'}, \, i, i'=1, 2, 3,\,  i\neq i'.$$
Similar to $r=2$, one may assume that 
$$0\neq \eta_1=m_1j_1+h_1j_2,$$
$$0\neq \eta_2=m_2j_1+h_2j_2,$$
$$0\neq \eta_3=m_3j_1+h_3j_2.$$
If the system of the corresponding $3$ equations of the form (**) is satisfied, then 
$$C_1(\eta_1^2-(\nu_1\cdot \omega^{(0)})^2)+C_2(\eta_2^2-(\nu_2\cdot\omega^{(0)})^2)+(\eta_3^2-(\nu_3\cdot\omega^{(0)})^2)=0,$$
where ($\bigstar$) is used to deduce that
$$C_1\eta_1+C_2\eta_2+\eta_3\equiv 0  \Leftrightarrow C_1\nu_1+C_2\nu_2+\nu_3=0,$$
cf. Lemma 2.1 in sect. 2.1 of \cite{W2}.

View $C_1$, $C_2$ as the unknown and 
write in the $j_1$, $j_2$ basis. In order for the $3$ equations to be compatible, the determinant $\Cal D$ of the $3\times 3$ matrix must satisfy:
$$\Cal D=\det \pmatrix m_1&h_1&\eta_1^2-(\nu_1\cdot \omega^{(0)})^2\\
m_2&h_2&\eta_2^2-(\nu_2\cdot \omega^{(0)})^2\\
m_3&h_3&\eta_3^2-(\nu_3\cdot \omega^{(0)})^2\endpmatrix=0.$$
Below we analyze the variety (in $m_i$, $h_i$, $i=1, 2, 3$) defined by $\Cal D=0$.  
 
 We first consider the case $$m_i\neq 0,\, h_i\neq 0,$$
 for all $i=1, 2, 3$. Let 
 $$\aligned\tilde \eta_i&=\frac{1}{m_ih_i}\eta_i,\\
 R_i&=\frac{1}{m_ih_i}(\eta_i^2-(\nu_i\cdot \omega^{(0)})^2))\\
 &=-(\frac{m_i}{h_i}+\frac{h_i}{m_i})+2(j_1\cdot j_2-\sqrt{j_1^2+1}\sqrt{j_2^2+1}), \endaligned$$
 for $i=1, 2, 3$. Without loss of generality, there are two cases:
 \item{a)} $\tilde\eta_2-\tilde\eta_1$, $\tilde\eta_3-\tilde\eta_1$ span $\Bbb R^2$; 
 \item{b)} $\tilde\eta_2-\tilde\eta_1$, $\tilde\eta_3-\tilde\eta_1$ are co-linear.
 
 \noindent Case a) Using row reduction, the $3\times 3$ determinant  (in the $j_1$, $j_2$ basis as before)
 $$\aligned D&=\det \pmatrix \tilde \eta_1&R_1\\
 \tilde\eta_2&R_2\\
 \tilde \eta_3&R_3\endpmatrix\\
 &=P_1( j_1\cdot j_2)-P_1 \sqrt {j_1^2+1}\sqrt{j_2^2+1}+P_2,\endaligned$$
 where $P_1$ and $P_2$ are rational functions of $m_i$, $h_i$, $i=1, 2, 3$, and 
 $$P_1=2\det (\tilde\eta_2-\tilde\eta_1, \tilde\eta_3-\tilde\eta_1)\neq 0.$$
 So $D\neq 0$, leading to a contradiction. 

\noindent Case b) Since  $$\tilde\eta_3-\tilde\eta_1= C(\tilde\eta_2-\tilde\eta_1), \, C\neq 0,$$
we have 
$$\frac{ \frac{1}{h_3}-\frac{1}{h_1} }{\frac{1}{h_2}-\frac{1}{h_1} }=\frac{ \frac{1}{m_3}-\frac{1}{m_1} }{\frac{1}{m_2}-\frac{1}{m_1} }=C,$$
assuming $m_1\neq m_2$ and $h_1 \neq h_2$.
Since $$\aligned D&=P_2\\
&=\big [-\big( \frac {m_2} {h_2} +\frac {h_2} {m_2}  \big)+ \big( \frac{m_1}{h_1}+\frac {h_1} {m_1} \big)\big] \det ( \tilde\eta_1, \tilde\eta_3-\tilde\eta_1)  \\
&\quad+ \big[\big( \frac {m_3} {h_3} +\frac {h_3} {m_3}  \big)- \big( \frac{m_1}{h_1}+\frac {h_1} {m_1} \big)\big] \det(\tilde\eta_1, \tilde\eta_2-\tilde\eta_1),\endaligned$$
if $D=0$, there is the additional equality
$$\frac{\big( \frac {m_3} {h_3} +\frac {h_3} {m_3}  \big)- \big( \frac{m_1}{h_1}+\frac {h_1} {m_1} \big)  }
{ \big(\frac {m_2} {h_2} +\frac {h_2} {m_2}  \big)- \big( \frac{m_1}{h_1}+\frac {h_1} {m_1} \big)}=C, \, C\neq 0,$$
assuming $\frac{m_1}{h_1}+\frac {h_1} {m_1} \neq \frac {m_2} {h_2} +\frac {h_2} {m_2}$.

We view $(m_1, h_1)$, $(m_2, h_2)$ as the given and write 
$$x=\frac{1}{m_3},\, y=\frac{1}{h_3}$$
as the unknown. We therefore have two equations of the form
$$\cases\frac{x-a_1}{a_2}=\frac{y-b_1}{b_2},\\
\frac {x-a_1} {a_2} =\frac{ (\frac{x}{y}+\frac{y}{x})-c_1 } {c_2},\endcases$$
where $a_i\neq 0$, $b_i\neq 0$, $c_i\neq 0$, $i=1, 2$ are the given.
Solving for $y$ from the first equation and substituting into the second
lead to a (non-zero) polynomial equation in $x\in \Bbb R$ of degree $3$. So there are at most $3$ 
solutions. Likewise, if $b_2=0$, then $a_2\neq 0$ and vice versa, 
and there are at most 2 solutions using similar arguments.

We now consider the remaining cases. Write $$\Cal D=\Cal P_1(j_1\cdot j_2)-\Cal P_1\sqrt{j^2_1+1}\sqrt{j_2^2+1}+\Cal P_2,$$
where $\Cal P_1$, $\Cal P_2\in \Bbb Z$, are functions of $m_i$, $h_i$, $i=1, 2, 3$.   
If $m_1h_1=0$ (without loss of generality, one may assume $m_1=0$), 
and $m_2h_2\neq 0$, $m_3h_3\neq 0$,
$\Cal P_1=0$ if and only if 
$h_3=h_2$. In that case 
$\Cal D=\Cal P_2=0$ has at most $1$ solution in $m_3$.
Finally if $m_1h_1=m_2h_2=0$, then $\Cal D\neq 0$. 

So in conclusion at most $5$ $\eta_i$'s of the form 
$$\eta_i=m_ij_1+h_ij_2$$
could possibly lead to compatible equations. Therefore
$$|S|=|\{(\nu_i, \eta_i)\}|\leq 5$$
and $$|\alpha|\leq 6.$$ 

Lastly, assume that there is a connected set $\alpha$ on $\Cal C_+$ with 
 $|\alpha|=2b+1$. Then there are $2b$ equations of the form (**). From the restrictions
 on the variable dependence of $\eta$, one can always choose a root so that $\eta$ in the 
 set $S$ contains a subset of $b+2$ elements of the forms 
 $$\aligned0\neq &\eta_1=m'_1j_1,\\
0\neq &\eta_2=m_1j_1+h_2j_2,\\
&\qquad\vdots\\
0\neq &\eta_\ell=m_1j_1+h_\ell j_\ell,\\
&\qquad\vdots\\
0\neq &\eta_b=m_1j_1+h_bj_b,\\
0\neq &\eta_{b+1}=m_1j_1+h_{b+1}j_\ell,\\
0\neq &\eta_{b+2}=m_1j_1+h_{b+2}j_\ell,
\endaligned$$
for some $\ell\in\{2, ..., b\}$. 
The arguments for $r=2, 3$ systems (when $r=3$, the ``remaining cases") show that there is no solution to the 
subsystem formed by $\eta_1$, $\eta_\ell$, $\eta_{b+1}$, $\eta_{b+2}$.
Therefore 
$$|\alpha|\leq 2(b-1)+1+1=2b.$$
Same arguments lead to 
$$|\alpha|\leq 2b,$$
on $\Cal C_-$.

To conclude, replace $p$ by $2p$ in (2.6) and define $\Gamma'$ to be 
$$\Gamma'=\text{supp }[(u^{(0)})^{*2p}]\backslash \{(0, 0)\}=\{(\nu, \eta)\}\subset\Bbb Z^{b+d}.\tag 2.6'$$ 
The above bounds on $|\alpha|$, on $\Cal C_+$ and $\Cal C_-$, clearly remain valid after replacing the $\Gamma$
in (2.6) by $\Gamma'$. We observe that if two points on $\Cal C_+$ are connected to the same point on $\Cal C_-$ by the $\Gamma$ in (2.6),
then they are connected on $\Cal C_+$ by the $\Gamma'$ in (2.6'); likewise after reversing the role of $\Cal C_+$ with $\Cal C_ -$.
Multiplying the bounds on $|\alpha|$ by $2$, therefore, proves the assertion.
\hfill$\square$
\enddemo

\noindent{\it Remark 1.} It is worth noting that the integer nature of $j_k$, $k=1, 2, ..., b$,
(aside from the non-degeneracy condition (ii)) is used in the variable reductions. This is 
contrary to the proof of the genericity of condition (iii) in the Lemma (sect.~2.1) and also 
that in sect.~2.2 of \cite{W2}.  

\noindent{\it Remark 2.} Clearly specializing to $\Theta=0$, this proves a weaker version 
of Proposition~2.1 with the bound $4b$ instead of $2d$. The geometric non-degeneracy condition (iii) is what 
permits the latter, in general, sharper bound on the characteristic $\Cal C$ defined in (1.9).  

To conclude, let $A$ be as in (2.5), following Lemma 5.2, there is the block decomposition for {\it every} $\Theta$ in (5.5): 
$$A_{\Cal C(\Theta)}=\oplus_{\alpha} A_{\alpha},\tag$\bigstar\bigstar\bigstar$ $$
where $\Cal C(\Theta)$ is the characteristic defined in ($\natural$), $\alpha$ are connected sets on $\Cal C(\Theta)$, and $A_\alpha$ are $A$ restricted to $\alpha$ -- therefore matrices of size at most $4b\times 4b$.

\subheading{5.4 Proof of Lemma 5.1}
\demo {Proof of Lemma 5.1}  Let $$\Theta\in\Cal Z,$$
the set defined in (5.5).
Since $|n|\leq N=|\log\delta|^s$ $(s>1)$, $\omega=\omega^{(0)}+\Cal O(\delta^p)$, $\Delta u^{(1)}=\Cal O(\delta^p)$ and $A(u^{(1)})=\Cal O(\delta^p)$, in view of (5.1, 5.5) and (D2), 
for small $\delta$, 
it suffices to look at $\theta$ such that 
$$\aligned \theta\in &\bigcup_{\Theta\in\Cal Z}  \{\Theta+\delta^{p}[-C |\log\delta|^{s}, C |\log\delta|^{s}]\}\\
:= &\bigcup_{\Theta\in\Cal Z} \Theta+I\\
:= &\,\Cal I,\endaligned$$
for some $C>1$. This is because, otherwise $T_N(\theta)$ is invertible, satisfying 
$$\Vert T_N^{-1}(\theta)\Vert\leq \delta^{-p}.$$

This can be seen as follows. From (5.1), for each $(n, j)\in [-N, N]^{b+d}$, on the diagonal, there is the quadratic polynomial in $\theta$: 
$$\Cal D'_{n, j} (\theta)= -(n\cdot\omega+\theta)^2+j^2+1,$$
where $\omega$ is the modulated frequency. The zeroes of $\Cal D'$:
$$\Theta'=-n\cdot\omega\pm\sqrt{j^2+1}.$$
The first variation at $\Theta'$:
$$\frac{\partial \Cal D'_{n, j} }{\partial \theta} (\Theta')=-2(n\cdot\omega+\Theta')=\pm 2\sqrt{j^2+1}.$$
So $$2\leq |\frac{\partial \Cal D'_{n, j} }{\partial \theta} (\Theta')|\leq 2(b+d)N=\Cal O( |\log\delta|^s),$$
for all $\Theta'$.
Since $|\omega-\omega^{(0)}|=\Cal O(\delta^p)$, 
the above, perturbation about the diagonals and (D2) yield the stated bound on 
$\Vert T_N^{-1}(\theta)\Vert$, for $\theta\not\in\Cal I$.

We now localize $\theta$ to the set $\Cal I$. 
Fix a $\Theta\in\Cal Z$ and write $\theta=\Theta+\delta^{p}\theta'$. 
Let $$w=\delta^{-1}a$$ and write $$\omega=\omega^{(0)}+\delta^p\omega'(w).$$ 
We have 
$$T_N(\theta)=\text{diag }[-[n\cdot \omega^{(0)}+\Theta+\delta^p(\theta'+n\cdot\omega')]^2+j^2+1]+\delta^{p} A_N (u^{(0)})
+\Cal O(\delta^{p+1}),$$
where $A_N$ is the restricted $A$ as defined in (2.5).

Let $\Cal Z_+$ be the set of $\Theta$ defined in (5.5) with the ``$+$" sign, and $\Cal Z_-$ the ``$-$" sign; 
$$\Cal Z_+\cup \Cal Z_-=\Cal Z.$$
For $\Theta\in\Cal Z_{\pm}$, define  
$$\Cal K= \text {diag }[\mp 2 \sqrt{j^2+1}(n\cdot\omega'+\theta')]+A_N.$$
Let $P$ be the projection onto $\Cal C(\Theta)$, and $P^c=I-P$.  Then expanding the diagonal and using that 
$$\delta ^{2p} |(n\cdot\omega'+\theta')^2| \leq \delta^{2p-1}$$
for $|n|$, $|\theta'|\leq \Cal O(|\log \delta|^s)$ and $\delta$ small enough, we arrive at 
$$P T_N (\theta)P= \delta^p P \Cal K P+ \Cal O(\delta^{p+1}).$$
It follows from Lemma 5.2, ($\bigstar\bigstar\bigstar$), that 
$$P\Cal K P=\oplus_\alpha\Cal K_\alpha(\theta', w),\tag 5.6$$
where $\theta'$ is in the size $\Cal O(|\log\delta|^s)$ interval $I$ introduced earlier, 
and $\Cal K_\alpha$ is $\Cal K$ restricted to $\alpha$, and therefore matrices of sizes at most $4b\times 4b$ for all $\alpha$.  

We proceed using the Schur reduction as 
in the proof of Lemmas~4.1. It suffices to estimate
$[PT_NP]^{-1}$,
as $P^cT_NP^c$ is invertible using (D2), and $$\Vert (P^cT_NP^c-\lambda)^{-1}\Vert\leq 4|\log\delta|^{s\Cal L}$$ uniformly in $\theta'$ for 
$$\lambda\in\big [-\frac{1}{4|\log\delta|^{s\Cal L}}, \frac{1}{4|\log\delta|^{s\Cal L}} \big].$$
Since $T_N=\Cal K+\Cal O(\delta^{p+1})$, this entails estimating 
$$[P\Cal K P]^{-1}=\oplus_\alpha[\Cal K_\alpha(\theta', w)]^{-1}.$$
The determinant of a block matrix $\det\Cal K_\alpha$ is a polynomial in $\theta'$ of degree at most 
$$M\leq 4b$$ with the coefficient in front of the highest degree term 
$C_M$ satisfying 
$$|C_M|>1.$$
Consequently, variation in $\theta'$ and summing over the number of $\Theta$ in $\Cal Z$, the set of zeroes defined in (5.5), satisfying $|\Theta|\leq \Cal O(|\log\delta|^s)$
proves (5.2, 5.4). Afterwards the point-wise estimates in (5.3)  follows. This is as in the proof of Lemma~4.1. 
\hfill $\square$ 
\enddemo

\item{{\bf (Fi)}} We now fix $\epsilon\in (0, 1/2)$ and $\epsilon'\in (0, 1)$ in Proposition~4.3;
subsequently $\sigma\in (0, 1) $ as well.

Before we establish the analogue of Lemma~5.1 for all scales, we first give a proof of Diophantine $\omega$, under 
conditions which will be verified along the iteration process in sect.~6.  

\subheading{5.5 Diophantine $\omega$}
In Proposition 4.3, (4.23), $\omega^{(1)}$ is shown to be Diophantine when restricted to the scale $N=|\log\delta|^s, s>1$,
using that $\omega^{(0)}$ is Diophantine and small $\delta$. Below we give a general proof 
for unrestricted $N$. 

\proclaim{Lemma 5.3}  Assume that $\omega^{(0)}$ is a Diophantine vector in $\Bbb R^b$ satisfying
$$\Vert n\cdot \omega^{(0)}\Vert_{\Bbb T}\geq\frac{2\xi}{|n|^\gamma},\quad n\in\Bbb Z^b\backslash\{0\},\,\xi>0, \gamma>2b.$$
Let $$\omega=\omega^{(0)}+\delta^p \omega',$$ where $\Cal O(1)=\vert\omega'\vert\leq 1$. There exists $\delta_0\in (0, 1)$, such that 
for all $\delta\in (0, \delta_0)$, 
$\omega$ is Diophantine satisfying
$$\Vert n\cdot \omega\Vert_{\Bbb T}\geq\frac{\xi}{|n|^{3\gamma}},\quad n\in\Bbb Z^b\backslash\{0\},\,\xi>0, \gamma>2b,\tag 5.7$$
away from a set in $\omega'$ of measure less than $\Cal O (\delta^p)$, where the order $\Cal O$ depends on $\xi$ and $\gamma$. 
\endproclaim
\demo{Proof}
If $$|n|\delta^{p}\leq\frac{\xi}{2|n|^\gamma},$$ clearly (5.7) holds. If 
$$|n|\delta^{p}>\frac{\xi}{2|n|^\gamma},$$ 
then 
$$|\delta^p n\cdot \omega'+n\cdot \omega^{(0)}+j|<\frac{\xi}{|n|^{3\gamma}}$$
for all $j\in\Bbb Z$, leads to a set in $\omega'$ of measure less than $\Cal O(\delta^{p})$.
This is because for each given $n$, $ j$ may be restricted to $|j|\sim |n|$. Summing over $(n, j)\in\Bbb Z^{b+1}$ then gives the 
measure estimate.
\hfill $\square$
\enddemo

\noindent{\it Remark.} Using Lemma 5.3 and Proposition 3.1, $\omega^{(1)}$ satisfies the unrestricted 
Diophantine property (5.7) with a further excision in $a$. 
\item{{\bf (Fii)}} We now fix $\xi$ and $\gamma$. 

\subheading{5.6 The iterative $\theta$ estimates}
We now establish the analogues of Lemma 5.1 for all scales $N$. 
Assume that for all $r\geq 0$,
$$|u^{(r)}(x)|<\delta e^{-\gamma |x|},\, \gamma>0.\tag 5.8$$
Let 
$$T_N(\theta)=T_N(\theta; u^{(r')}), \text { with } r'=C'\log N,\, C'>1.\tag 5.9$$
Assume that
$$\Vert u^{(R)}-u^{(r')}\Vert <\delta^pe^{-\tilde\gamma N},\, \tilde\gamma>\gamma>0,\tag 5.10$$
for all $R\geq r'$. 

\proclaim{Proposition 5.4}
Let $u^{(0)}=\sum_{k=1}^b a_k \cos ( {-(\sqrt{j_k^2+1})t}+j_k\cdot x)$
be a solution to the linear Klein-Gordon equation (1.2),  satisfying the non-degeneracy conditions (i, ii), and (4.17, 4.19) hold
with $\epsilon\in (0, 1/2)$.
Let $\sigma, \tau$ be numerical constants satisfying 
$$0<\tau<1/s<\sigma<1.\tag $\spadesuit$ $$ 
There exist $\delta_0>0$ and $0<\beta<\gamma<\tilde\gamma$, such that for all $\delta\in(0, \delta_0)$, all $N>10/\beta^2$,
$$\Vert [T_N(\theta)]^{-1}\Vert<e^{N^\sigma}, \tag 5.11 $$
$$|[T_N(\theta)]^{-1}(x,y)|<e^{-\beta |x-y|},\tag 5.12$$
for all $x, y$ such that $|x-y|>1/\beta^2$, for $\theta$ away from a set $B_N(\theta)\subset\Bbb R$ with 
$$\text{meas }B_N(\theta)<e^{-N^\tau}.\tag 5.13$$
\endproclaim 
\demo{Proof} For scales $N\in (10/\beta^2, e^{|\log\delta|^{3/4}})$, one may clearly proceed as in the proof of Lemma~5.1.

For $N\geq e^{|\log\delta|^{3/4}}$, pave $[-N, N]^{b+d}$ by 
translates of the cube $[-N^{1/C'}, N^{1/C'}]^{b+d}$, 
with $1<C'<1/\tau$. The key is to control the 
number of {\it bad} $N^{1/C'}$-boxes in $[-N, N]^{b+d}$
for (any) fixed $\theta$, on which (5.11) or (5.12) is violated. 
According to the established mechanism such as that in \cite{BGS},
this needs to be at most sublinear in $N$ for an iterative proof of (5.11) and (5.12)
at all scales $N$. The proof is similar to that of Lemma~4.1 in \cite{BW}. We give the gist below.

Define $$T=[-N, N]^b\times [-N^{1/C'}, N^{1/C'}]^d.$$ 
Requiring
$$|(n\cdot\omega+\theta)^2-j^2-1|>e^{-N^{1/2C'}},$$
for all $(n, j)\in [-N, N]^{b+d}\backslash T$, excises
a set $\Theta$ in $\theta$  of measure satisfying 
$$\text {meas }\Theta\leq N^{C(b+d)} e^{-N^{1/2C'}}\leq e^{-N^{1/3C'}}.$$
Using the decaying nonlinear term, diagonal perturbation and Neumann series give that 
for $\theta\notin \Theta$, (5.11) and (5.12) are satisfied, on all $N^{1/C'}$-boxes $\Lambda\subset  [-N, N]^{b+d}\backslash T$, if $1/2C'>\sigma$. 
In particular, there are {\it no} bad $\Lambda$'s for any fixed $\theta\notin\Theta$.

For $\Lambda$, such that $\Lambda\cap T\neq\emptyset$, 
Diophantine $\omega$ and the $\theta$ estimate (5.13), yield that for any fixed $\theta$, there
are at most sublinear in $N$ such $\Lambda$, on which (5.11) or (5.12) are violated. 
This uses the covariance structure: 
$$\theta\mapsto n\cdot\omega,$$
and can be seen as follows. Assume that (5.11) and (5.12) hold at scale $N^{1/C'}$. 
The former maybe written as an algebraic inequality in $\theta$, by using the Hilbert-Schmidt 
norm; the latter is the division of two determinants, and hence algebraic. These algebraic 
inequalities are of degree at most $N^{9(b+d)/C'}$. So $B_{N^{1/C'}}:=\Cal B$ is semi-algebraic and 
has at most $N^{10(b+d)/C'}$ connected components, by using Basu's theorem \cite{Ba}, restated 
as Theorem~7.3 in \cite{BGS}. 

If there are $n$ and $n'$, $n\neq n'$, such that $n$ and $n'$ belong to the
same connected component of $\Cal B$, then 
$$|(n-n')\cdot\omega|<e^{-N^{\tau/C'}},$$
from (5.13) and using the covariance. But this contradicts the Diophantine estimate:
$$|(n-n')\cdot\omega|\geq \frac{\xi} {N^{3\gamma}},$$
from (5.7). Therefore there can be at most {\it one} integral points in each connected component,
which leads to at most sublinear in $N$ bad $N^{1/C'}$-boxes. (Cf. the proof of Lemma~4.1 \cite{BW}.)

Combining the two regions, gives that for any fixed $\theta\notin \Theta$, there are only sublinear $\Lambda$,
on which (5.11) and (5.12) are violated. Standard arguments (by now), such as those debuted in \cite{BGS}, conclude the proof by choosing 
$C'<1/3\tau$.
\hfill $\square$
\item{{\bf (Fiii)}} We fix $\tau$ and $s$ satisfying ($\spadesuit$). 
\enddemo

\head {\bf 6. Proof of the Theorem}\endhead
Using a Newton iteration to solve the $Q$ and $P$-equations, the proof of the Theorem is an induction.
It is based on the initial corrections in Proposition 4.3, the iterative $\theta$ estimates in Proposition~5.4, the covariance structure
$$\theta\mapsto n\cdot\omega,$$
the amplitude-frequency map: $\omega=\omega(a)$, resolvent expansions and additional excisions in the amplitude $a$. 
This is essentially the same as sects. 5 and 6 \cite{BW}, which deals with a decaying nonlinear term.

Let $$w=\delta^{-1}a\in (-1, 1)^b\backslash\{0\},\, \omega\in \delta^{p}(-K, K)^b+(\sqrt{j^2_1+1},\sqrt{ j^2_2+1}, ..., \sqrt{j^2_b+1}),$$
where $K=K(p, b)$ and define
$$\tilde\omega=[\omega- (\sqrt{j^2_1+1},\sqrt{ j^2_2+1}, ..., \sqrt{j^2_b+1})/\delta^{p}\in (-K, K)^b.$$
Let $M$ be a large integer.
The proof consists in showing that on 
the {\it entire} $(w,\tilde\omega)$ space, namely $(-1, 1)^b\backslash\{0\}\times (-K, K)^b$,
the following
assumptions are verified for all $r\geq 1$ and fixed $\delta$ sufficiently small: 
\item{(Hi)} $\text{supp }u^{(r)}\subseteq B(0, M^r)$ ($\text{supp }u^{(0)}\subset B(0, M)$)
\item{(Hii)} $\Vert \Delta u^{(r)}\Vert<\delta_r$, $\Vert \partial \Delta u^{(r)}\Vert<\bar\delta_r$, with $\delta_{r+1}\ll\delta_{r}$ and $\bar\delta_{r+1}\ll\bar\delta_{r},$
where $\partial$ refers to derivations in $w$ or $\tilde\omega$ and 
$\Vert\,\Vert:=\sup_{w,\tilde\omega}\Vert\,\Vert_{\ell^2(\Bbb Z^{b+d})}$. (See the precise bounds at the conclusion of the proof.)
\item{(Hiii)} $|u^{(r)}(x)|<\delta e^{-\alpha |x|}$ for some $\alpha>0$
\smallskip
Using (Hi-iii), an application of the implicit function theorem to the
$Q$-equations: 
$$\tilde\omega_k(w)=\frac{[(\Cal M{u})^{*p+1}\Cal M ] (-e_k, j_k)}{2^{p+1} w_k \sqrt{j_k^2+1} },\tag 6.1$$
$k=1, 2, ..., b$, with $u=u^{(r)}$, a {\it real} valued $C^1$ function, yields
$$ {\tilde\omega}^{(r)}_k(w)=\Omega_k(u^{(0)}(w))+\delta^{p-1} \phi_k^{(r)}(w)\tag 6.2$$
where the polynomials $\Omega_k$, $k=1,2, ..., b$ are as in (3.1), homogeneous in $w$ of degree $p$, $0<\epsilon<1/2$ in view of (4.18), and  
$\Vert\partial\phi^{(r)}\Vert<C$. 

We define $\phi_0=0$ and denote the graph of $ {\tilde\omega}^{(r)}$ by $\Phi_r$. 
Moreover  by (Hii),
$$| {\tilde\omega}^{(r)}-  {\tilde\omega}^{(r-1)} |\leq \Cal O(1)\Vert u^{(r)}-u^{(r-1)}\Vert <\delta_r,\tag 6.3$$
so that $\Phi_{r-1}$ is an $\delta_r$ approximation of $\Phi_{r}$. (Cf. the proof below (5.10) in \cite{W2}.) \hfill $\square$
\smallskip
Below we continue with the assumptions on the {\it restricted} intervals in $(w,\tilde\omega)$ on $(-1, 1)^b\backslash\{0\}\times (-K, K)^b$,
where approximate solutions could be constructed.
\item{(Hiv)} There is a collection $\Lambda_r$ of intervals of size $cM^{-r^C}\delta^\epsilon$, $\epsilon\in (0, 1/2)$, $C>1$ such that 
\item{(a)} On $I\in\Lambda_r$, $u^{(r)}(w,\tilde\omega)$ is given by a rational function in $(w,\tilde\omega)$ of degree at most 
$M^{Cr^3}$, 
\item{(b)} For $(w,\tilde\omega)\in\bigcup_{I\in\Lambda_r} I$,

$\Vert F(u^{(r)})\Vert<\kappa_r$,
$\Vert \partial F(u^{(r)})\Vert<\bar\kappa_r$ 
with $\kappa_{r+1}\ll\kappa_{r}$ and $\bar\kappa_{r+1}\ll\bar\kappa_{r}$
\item{(c)} Let $N=M^r$. For $(w,\tilde\omega)\in\bigcup_{I\in\Lambda_r} I$, $T=T(u^{(r-1)}):=F'(u^{(r-1)})$ satisfies

$\Vert T_N^{-1}\Vert <M^{(r^C+|\log\delta|})$, 

$|T_N^{-1}(x, y)|<e^{-\beta |x-y|}, \, \beta>0$, for $|x-y|>Cr^{C}$,

where $T_N$ is $T$ restricted to $[-N, N]^{b+d}$.
\item{(d)} Each $I\in\Lambda_r$ is contained in an interval $I'\in\Lambda_{r-1}$ and 
$$\text{meas}_b (\Phi_r\cap(\bigcup_{I'\in\Lambda_{r-1}} I'\backslash \bigcup_{I\in\Lambda_{r}}I)<\delta^{c\epsilon}[\exp\exp(\log (r+1))^{1/3}]^{-1},\, r\geq 2.$$
If $\tilde\omega\in \Phi_r\cap I$, then $$\omega=\omega^{(0)}+\delta^p\tilde\omega$$
is Diophantine satisfying
$$\Vert n\cdot\omega\Vert_{\Bbb T}\geq \frac{\xi}{|n|^\gamma},\, \xi>0,\, \gamma>6b,$$
for $|n|\leq M^r$,
after identification of $\Phi_r\cap I$ with an interval in $\Bbb R^b$. 

We remark that the approximate solutions $u^{(r)}$ are defined, a priori, on $\Lambda_r$, but as $C^1$ {\it functions} they can be 
extended to $(-1, 1)^b\times (-K, K)^b$, using a standard extension argument, cf. sect. 10, (10.33-10.37) in \cite{B3},
thus verifying (Hi-iii) (and hence (5.8)). This is also why the amplitude-frequency map $a\mapsto \omega(a)$ is a diffeomorphism on $(-\delta, \delta)^b$
in the Theorem.
\smallskip
\noindent {\it Proof of the Theorem}

As mentioned earlier, the induction follows that in sects.~5 and 6 \cite{BW}, below we emphasize a few key steps.

The first $R$, $R=|\log\delta|^{3/4}$, 
steps of the induction are provided by the generalization of Proposition 4.3 to include all 
scales $$N\in [M,  M^{|\log\delta|^{3/4}}].$$
The iteration to subsequent scales
uses Proposition~5.4.  

Let $u$ denote $u^{(0)}$, $u^{(1)}$, ... For all $\bar N$, let $T_{\bar N}=T_{\bar N}(u)$ be the linearized operator 
evaluated at $u$ and restricted to $\{J+[-\bar N, \bar N]^{b+d}\}\times \{0, 1\}$, where $J\in[-N, N]^d$. (For simplicity the 
$J$ subindex is omitted.) Define the operator $T_{\bar N}(\theta)$ as before. Assume that (Hi-iv) hold at stage $r$. 
On the set of intervals $\Lambda_r$ in (Hiv), there are moreover the following estimates
from Proposition~5.4.

\proclaim{Lemma 6.1} 
$$\align\Vert T_{\bar N}^{-1}(\theta)\Vert &<e^{{\bar N}^\sigma},\tag 6.4\\
|T_{\bar N}^{-1}(\theta)(x, y)| &< e^{-\beta |x-y|},\, \beta>0,\tag 6.5\endalign$$
for all $x, y$ such that $|x-y|>\bar N/10$, away from a set $B_{\bar N}(\theta)$ with 
$$\text{meas } B_{\bar N}(\theta)<e^{-\bar N^\tau},$$
where $u=u^{(r)}$, $|\log\delta|^s\leq \bar N\leq r^C$, $C>2/\sigma$, 
$r\geq R$.
\endproclaim
Assume that (Hi-iv) hold at step $r$, the iteration to step $r+1$ paves the cube $[-M^{r+1}, M^{r+1}]^{b+d}$ by the cube
$[-M^{r}, M^{r}]^{b+d}$ and much smaller cubes $\Lambda$, of size 
$$M_0=(\log M^r)^C\sim r^C.$$ Define 
$$T=[-M^{r+1}, M^{r+1}]^b\times [-M_0, M_0]^d.$$ In the paving process, there are two types of $\Lambda$:
\item{(i)} $\Lambda\cap T\neq\emptyset$;
\item{(ii)} $\Lambda\cap T=\emptyset$.

Region (i) uses Lemma~6.1 and semi-algebraic projection, Lemma~9.9 \cite {B4}. After removing a set in $a$ of measure at most
$M^{-r/10}$, (6.4) and (6.5) hold for all such $\Lambda$. Concretely, we first divide the $a$ parameter space into sufficiently small intervals. On each such interval, we make rational approximations to $a=a(\omega)$,
using the $Q$-equations at $u=u^{(0)}$. 
Using this truncated $a$, the Newton iteration, which uses resolvent series,
preserves rationality, and leads to  
$u^{(r)}$ rational in $\omega$ for all $r$. Here we used also that
the $P$-equations only depend explicit on $\omega$ and $u$, and have {\it no explicit} dependence on $a$.
Afterwards, one may follow the arguments in 
Chap.~18, (18.28)-(18.33) \cite{B4}. 
Region (ii) makes direction excisions in $\omega$ ($a$). After removing a set of measure at most 
$$e^{-M_0^\sigma} M^{r+1}< e^{-r^{C\sigma}} M^{r+1}<e^{-r^2},$$
by choosing $C>2/\sigma$, (6.4) and (6.5) hold on all such $\Lambda$. Resolvent equation then yields (Hiv, c). 
(Cf. (5.8)-(5.20) \cite{BW}.)
Afterwards, we may reproduce (Hiv) at step $r+1$. Using the (Hiv, c) in the Newton iteration, then proves (Hi-iii) at
step $r+1$. This iteratively solves the $Q$
and $P$-equations. Moreover there are the bounds
$$\delta_r<\delta^{p} M^{-(\frac{4}{3})^r}, \, \bar\delta_r<\delta^{p} M^{-\frac{1}{2}(\frac{4}{3})^r}; \kappa_r<\delta^{3p-1} M^{-(\frac{4}{3})^{r+2}}, \, \bar\kappa_r<\delta^{3p-1} M^{-\frac{1}{2}(\frac{4}{3})^{r+2}},$$
$\epsilon\in (0, 1/2)$, cf. the proofs of Lemmas 5.5 and 5.2 in \cite{W2}. 

From Proposition 3.1, for small $\delta$, the (closed) set $\Cal B'$ has a semi-algebraic description, in terms of 
one (non-zero) polynomial in $w$ of degree $b(p-1)$, namely $\det (\partial \Omega/\partial w)$, where 
$\Omega=\{\Omega_k\}_{k=1}^b$ and $\Omega_k$ as in (3.1). (Recall that $w=\delta^{-1} a$.)
The number of connected components
of $\Cal B'$ is therefore at most $C_bb^b(p-1)^b$, cf. Theorem 1 in \cite{Ba} or Theorem 9.3 in \cite{B4}. 
The complement contains an open set of measure at least 
($1-\epsilon'/2$) in $w$. 

The proceeding iterative construction with the amplitude-frequency diffeomorphism:  $a\mapsto \omega(a)$
mentioned after (Hiv, d) then proves the Theorem. 
on a Cantor set of measure at least
($1-\tilde \epsilon$) for some $0<\epsilon'<\tilde \epsilon<1$. Renaming
$\tilde\epsilon$, $\epsilon$ concludes
the proof.  
\hfill$\square$

\enddocument
\end